\documentclass[11pt]{amsart}
\usepackage{amssymb}
\usepackage{graphicx}

\newcommand{\Z}{{\mathbb Z}}
\newcommand{\R}{{\mathbb R}}
\newcommand{\Hy}{{\mathbb H}}
\newcommand{\E}{{\mathbb E}}
\newcommand{\C}{{\mathbb C}}
\newcommand{\Diff}{\mathop{\text{Diff}}}
\newcommand{\goesto}{\longmapsto}

\newcommand{\df}[1]{\emph{#1}}
\newcommand{\boundary}{\mathop {\partial}}
\newcommand{\simplex}[1]{\ensuremath{ \left \langle #1 \right \rangle }}
\newcommand{\set}[1]{\ensuremath{\left \{ #1 \right \} }}

\newcommand{\degree}{^\circ}
\newcommand{\restrict}{\,\big\vert\,}
\newcommand{\floor}[1]{ \left \lfloor #1 \right \rfloor}
\newcommand{\ceiling}[1]{ \left \lceil #1 \right \rceil}

\DeclareMathOperator{\PSL}{\text{PSL}}
\DeclareMathOperator{\SL}{\text{SL}}
\DeclareMathOperator{\SO}{\text{SO}}
\DeclareMathOperator{\Homeo}{Homeo}
\DeclareMathOperator{\TS}{TS}
\DeclareMathOperator{\MG}{MG}
\DeclareMathOperator{\PG}{PG}
\DeclareMathOperator{\CP}{\mathbb{CP}}
\DeclareMathOperator{\RP}{\mathbb{RP}}
\DeclareMathOperator{\cl}{cl} 
\DeclareMathOperator{\acl}{acl} 
\DeclareMathOperator{\gfl}{Gfl} 
\newtheorem{theorem}{Theorem}[section]
\newtheorem{corollary}[theorem]{Corollary}
\newtheorem{lemma}[theorem]{Lemma}
\newtheorem{proposition}[theorem]{Proposition}
\newtheorem{conjecture}[theorem]{Conjecture}
\theoremstyle{definition}
\newtheorem{definition}[theorem]{Definition}
\theoremstyle{remark}
\newtheorem{remark}[theorem]{Remark}
\newtheorem{example}[theorem]{Example}

\begin{document}
\title[Three-manifolds, Foliations and Circles I]
{Three-manifolds, Foliations and Circles, I\\
Preliminary version}
\author{William P. Thurston}
\date{\today}
\address{
Mathematics Department \\
University of California at Davis \\
Davis, CA 95616
}
\email{wpt@math.ucdavis.edu}
\begin{abstract}
A manifold $M$ \df{slithers} around a manifold $N$ 
when the universal cover of $M$ fibers over $N$ so that deck
transformations are bundle automorphisms. 
Three-manifolds that slither around $S^1$ are like a hybrid
between three-manifolds that fiber over $S^1$ and certain
kinds of Seifert-fibered three-manifolds. There are
examples of non-Haken hyperbolic manifolds that slither around $S^1$. 
It seems conceivable that every hyperbolic 3-manifold
slithers around $S^1$, and it seems reasonable that
every hyperbolic three-manifold has a finite sheeted cover
that slithers around $S^1$.

If $M$ is a closed $3$-manifold, then
\textbf{I.}
$M$ slithers around the circle if and only if it 
has a \df{uniform} foliation $\mathcal F$, defined to be
a foliation without Reeb components such that 
in the universal cover any two leaves are a uniformly bounded distance
apart.

\textbf{II.}
Every uniform foliation $\mathcal F$
has a transverse flow $\phi_t$ 
that is either pseudo-Anosov, periodic, or reducible (admits a
non-empty collection of invariant incompressible tori and
Klein bottles).

\textbf{III.}
If $M$ is hyperbolic and $\mathcal F$ is a uniform foliation of $M$,
the stable and unstable laminations for $\phi_t$
are quasi-geodesic.  The leaves of $\mathcal F$ extend continuously to give
$\pi_1(M)$-equivariant sphere-filling curves
in the sphere at infinity of $\tilde M$.

\textbf{IV.} 
The \df{skew $\R$-covered Anosov foliations}
analyzed by S\'ergio Fenley \cite{Fenley:skew} 
slither around the circle.  They correspond 1--1 to cocompact
\df{extended convergence groups}, which are
subgroups $\Gamma \subset \widetilde {\Homeo(S^1)}$ such that
$\Tilde T/\Gamma$ is Hausdorff,  where
$T$ is the set of counter-clockwise ordered triples of distinct points
on the circle. (Convergence groups are the special case that $\Gamma$
contains the kernel $\Z \to \widetilde {\Homeo(S^1)} \to \Homeo(S^1)$.)

\textbf{Preview.}
Two or more further parts are projected in this series. Part II
will analyze  the asymptotic geometry of leaves of taut foliations of
3-manifolds and construct a universal circle-at-infinity that
collates the circles-at-infinity for all the leaves.
Provided that $M$ is atoroidal, the action of $\pi_1(M)$ on this circle
will be used to construct a genuine essential lamination transverse to any
taut foliation. 

In a subsequent part, we plan to prove the geometric
decomposition conjecture for three-manifolds that slither around $S^1$ by
analyzing the deformation theory of uniform `quasi-Fuchsian' 
foliations of $M\times \R$ whose  leaves have three-dimensional hyperbolic
structures.
\end{abstract}
\maketitle
\tableofcontents
\listoffigures
%

\section{Fiberings and Slitherings}

\begin{definition} \label{dfn: slither}
One manifold $M$ \df{slithers} around a second manifold $N$ when there  is
a fibration $s: \Tilde M \to N$  of some regular covering space
$p: \Tilde M \to M$ whose
deck transformations are bundle automorphisms for $s$. In other
words, deck transformation take
each fiber of $s$ to a (possibly different) fiber of $s$.  This
structure, determined by $s$, is a \df{slithering}.

The manifold $M$ is the \df{total space}, and $N$ is the 
\df{base}. The \df{fibers} of a slithering are the fibers of $s$.
The components of the
images $p(F)$ in $M$ of the fibers are the
\df{leaves} of the slithering, and they form a foliation $\mathcal F(s)$. 
\end{definition}

We could always use the universal cover of $M$ for the covering
space $\Tilde M$, but it is sometimes convenient to construct examples
in terms of other regular covering spaces.  A fibration $M^m \to N^n$
qualifies as a slithering, but there are many examples that are not of this
form. To start,
\begin{example} \label{example: torus slithering}
\begin{figure}[hbp] 
\centering
\includegraphics[width=.55\textwidth]{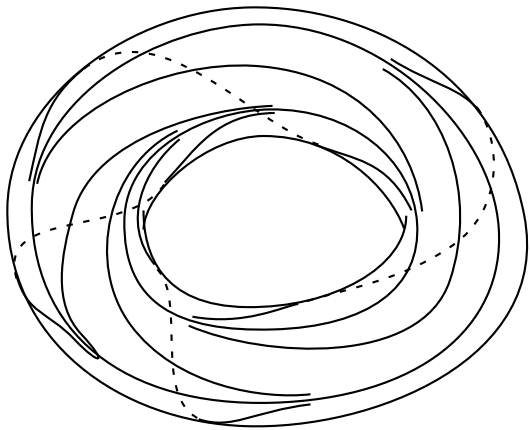}
\caption[A slithering of $T^2$]
{ \label{figure: torus slithering}
This foliation of $T^2$ has two closed leaves; all
other leaves are homeomorphic to $\R$, spiraling to the two
closed leaves at their two ends. Its universal cover
can be represented as the foliation of $\R^2$ by horizontal lines, where
the deck transformations are generated
(for example) by $(x,y) \goesto (x,y+2\pi)$ and
$(x,y) \goesto (x+1, y + .5 \sin(y))$.  These transformations
 act as automorphisms
of the fibration $\R^2 \to S^1 = \R^2 / (\R \oplus 2\pi\Z)$.
}
\end{figure}
The only closed 2-manifolds that can slither are the torus and Klein bottle,
which fiber over $S^1$.  However, these manifolds also have slitherings
that are not fiberings. For instance, figure \ref{figure: torus slithering}
shows a foliation with two closed leaves, where all other leaves are
lines spiraling to the two closed leaves in the two directions. The
universal cover of $T^2$ can be represented as $\R^2$ with
the foliation by horizontal lines, where deck transformations can
be taken as the group generated by 
\[
\phi(x,y) = (x+1, y+.5\sin(y))  \quad \psi(x,y) = (x,y+2\pi).
\]
This group acts as automorphisms of the fibering over $S^1$
\[
\R^2 \to \R^2/(\R \times  2 \pi \Z) = S^1.
\]

The particular fibering over $S^1$ is part of the data, and is not
determined by the foliation (although in many examples,
there is a unique simplest choice.)  One could, for example,
use the fibering over the $2 \pi k$ circle $\R^2/(\R\times 2 \pi k \Z)$.
\end{example}

\begin{figure}[thbp]
\begin{minipage}{.30\textwidth}
\includegraphics[width=\textwidth]{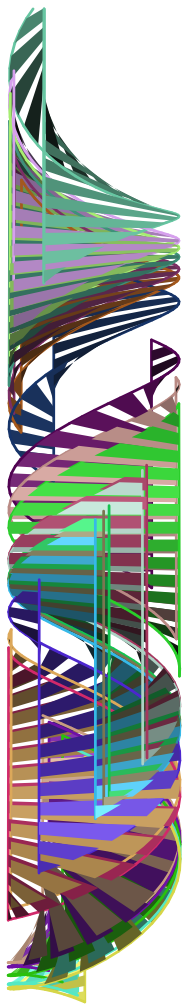}
\end{minipage}
\begin{minipage}{.7\textwidth}
\caption[The circle at infinity foliation of $TS(\Hy^2)$]
{\label{figure: circle at infinity}
This is the universal cover of the tangent
circle bundle of $\Hy^2$ (or any hyperbolic surface), with a few randomly selected
leaves of its circle-at-infinity foliation. Each leaf is rendered as an
assemblage of slats, so we can glimpse what's behind.
This foliation is also called the stable
foliation for the geodesic flow.  The horizontal direction
is the projective disk model for $\Hy^2$, while the vertical direction
is the angle of a vector in the plane (using ordinary $\E^2$ measurement.)
In these coordinates, one point at infinity on each leaf blows
up to a vertical interval.  The leaves are portions of helicoids
that the foliation has cleverly stacked.
Any projective automorphism
of the disk has a derivative that when lifted to $\Tilde TS(\Hy^2)$
preserves this foliation.
The unstable foliation for the geodesic
flow is obtained by rotating this picture $180\degree$ about its vertical axis.
}
\end{minipage}
\end{figure}
\begin{example} \label{example: tangent of hyperbolic}
Let $Q^n$ be a hyperbolic manifold or orbifold, and let $\TS(Q^n)$
be its tangent sphere bundle. Then  $\TS(\Hy^n) \to \TS(Q^n)$
is a regular covering, and the map $s: \TS(\Hy^n) \to S_\infty^{n-1}$
that sends each tangent ray to its endpoint at infinity
is a fibration. The deck transformations act as bundle automorphisms,
so $\TS(Q^n)$ slithers around $S_\infty^{n-1}$.
\begin{figure}
\begin{minipage}{.32\textwidth}
\includegraphics[width=\textwidth]{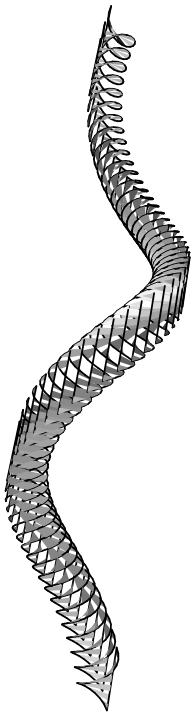}
\end{minipage}
\begin{minipage}{.65\textwidth}
\caption[The circle at infinity, condensed]{
\label{figure2: circle at infinity}
This is the same foliation of $\Tilde \TS(\Hy^2)$ shown
in figure \ref{figure: circle at infinity}, but
drawn with every leaf scaled $\times .25$ toward the point where its
geodesics are converging, in order to
bring out the helical relationship of leaves to points on a circle at infinity.
}
\end{minipage}
\end{figure}

When $n=2$, the restriction of the bundle $\TS(Q^2) \to Q^2$
to any geodesic in $Q^2$ is a torus, and
the slithering of $\TS(Q^2)$ around $S^1$
induces a slithering of this torus around $S^1$ that is
topologically equivalent to the first slithering of the preceding example.
\end{example}

\begin{figure}[htbp]
\begin{minipage}{.45\textwidth}
\includegraphics[width=\textwidth]{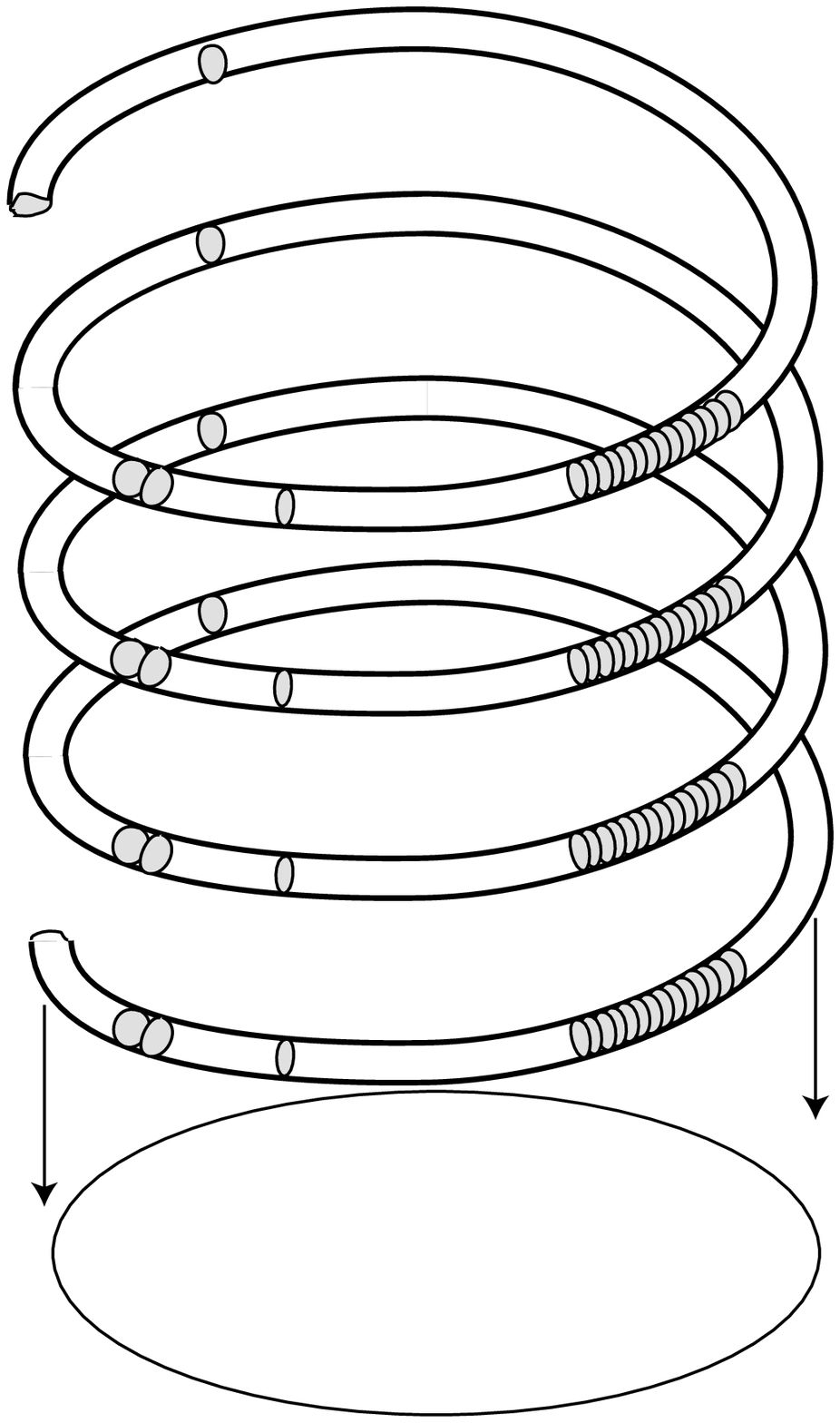}
\end{minipage}
\begin{minipage}{.55\textwidth}
\caption[Slitherings around $S^1$ give a long coiled shape to $\tilde M$.]
{\label{figure: snake}
A slithering around $S^1$ can be pictured as giving a long, coiled shape
to the universal cover $\Tilde M$. The fundamental domains for the
action of deck transformations have projections to
$S^1$ that are isomorphic, but for a `typical' slithering no two
fundamental domains are in phase.  The leaves of the foliation
of $\Tilde M$ group into families (a few of them shown) that
are the fibers of the slithering. Deck transformations of $\Tilde M$
move each family as a unit but alter the spacing between fibers.
Nonetheless, the distance between any
two leaves is bounded by a constant times the
number of turns of the coil (\ref{proposition: slithering uniform}.)
}
\end{minipage}
\end{figure}

When $M^m$ slithers around $N^n$, then the slithering lifts
to a slithering around $\Tilde N^n$.  In particular, if the
base is $S^1$, then the universal cover of
$M^m$ fibers over $\R$. One can picture a long stretched-out
image of $\Tilde M^m$, coiled around and around the circle (figure
\ref{figure: snake}.) The
fibers of the fibration to $S^1$ have infinitely many components.
Deck transformations of $\Tilde M^m \to M^m$
are periodic, but they probably do not move in
a uniform way: in some places the fibers squeeze closer together,
while elsewhere they spread apart, 
reminiscent of the slithering, undulating motion of a snake.
A slithering of a manifold
around $S^1$ has a strong dynamic element, stemming from the
hidden action of the fundamental group of $S^1$.
Slitherings are sneaky structures that you wouldn't be likely
to see if you weren't on watch for them.

A \df{foliated bundle} is a fibration $F^{m-p} \to M^m \to P^p$
with a reduction to a discrete structure group $\Homeo_\delta(F^{m-p})$,
that is, a foliation of dimension $p$ transverse to the fibers
whose leaves project as covering spaces to $P^p$.  Any such foliation,
when pulled back to the universal cover $\Tilde P^p$ of the base,
is a product $\Tilde P^p \times F^{m-p}$. In other words, the leaves
define a fibration over $F^{m-p}$, so $M^m $ slithers around
$F^{m-p}$.

The base can also be an orbifold, and the same reasoning applies. 
There is a well-developed theory started by Milnor (\cite{Milnor:inequality})
and Wood (\cite{Wood:inequality} concerning
which circle bundles admit foliations transverse to the fibers
when the base is a 2-dimensional orbifold (see
section \ref{section: groups, inequalities and topology}.)

Here is another way to represent the data for a slithering, in
a compact form that does not mention $\Tilde M^m$.  When $M^m$
slithers around $N^n$, then
$\pi_1(M^m)$ acts as a group of homeomorphisms of $N^n$.  There is a foliated
bundle $N^n \to E^{n+m} \to M^m$ associated with this action, obtained by
taking $\Tilde M^m \times N^n$ modulo the diagonal action.  The graph of the
fibration $\Tilde M^m \to N^n$ is invariant by the action of $\pi_1(M^m)$, so
it descends to give a section $M^m \to E^{n+m}$, transverse to the foliation of
$E^{n+m}$, and inducing the foliation of $M^m$.

For example, assuming $M^3$ is compact, a fibration $M^3 \to S^1$ is the same
thing as a section of the bundle $M^3 \times S^1 \to M^3$ that is transverse to
the horizontal foliation by $M^3 \times \theta$ associated with the trivial
action.  As a second example, when $M^3$ is a foliated circle bundle, the
fibration can be pulled back to the total space, giving a foliated circle
bundle together with a canonical section.

Every slithering gives data of this form, and the data is sufficient
to reconstruct the slithering.  What's often not obvious from this
type of data is whether or not the map $\Tilde M^m \to N^n$ is
actually a fibration; this depends on the global structure of
$\mathcal F(s)$.

For example,
a slithering over $\R$ gives a codimension one foliation such that
the space of leaves of the universal cover is homeomorphic to $\R$.
However, not every such foliation is a slithering over $\R$. For example,
consider $\R^3 \setminus \set 0$, modulo the action of $\Z$ generated
by $X \goesto 2X$. The quotient is $S^2 \times S^1$.
The foliation of $\R^3$ by horizontal planes
restricted to a foliation of the universal cover of $S^2 \times S^1$
such that the space of leaves
is $\R$; however, this map does not give a fibration over $\R$.
Furthermore, this foliation can be described as the foliation
induced from a section of a foliated circle bundle over $S^2 \times S^1$,
the bundle whose fiber is the one-point compactification of the
space $\R$ of leaves in the universal cover.

Among three-manifolds, the example of $S^2 \times S^1$ is exceptional.
Here is a fact from foliation theory:

\begin{proposition}\label{proposition: R-covered R-slitherings}
Let $\mathcal F$ be a codimension one foliation of an irreducible three-manifold
$M$. If the space of leaves of the foliation $\Tilde {\mathcal F}$
in the universal cover $\Tilde M$ is homeomorphic to $\R$, then
$\Tilde M$ is homeomorphic to $\R^3$, in a way that takes
the foliation $\Tilde {\mathcal F}$ to the foliation of $\R^3$ by horizontal
planes. In other words, $M$ slithers around $\R$.
\end{proposition}
Much more is actually known: Palmeira \cite{Palmeira:OpenFoliations}
showed that any foliation of an open $3$-manifold  by planes
is homeomorphic to the product of $\R$ with a foliation of the plane
(and similar results in higher dimensions). Poincar\'e studied 
foliations (and vector fields) in the plane, and showed that every leaf
of a foliation is a properly embedded line. It is easy to deduce that
if the space of leaves of a foliation of the plane is homeomorphic to
$\R$, then the leaves are fibers of a fibration.  
Haefliger classified all possible foliations of $\R^2$
in terms of the space of leaves, which is
a simply-connected but non-Hausdorff $1$-manifold, together with
with certain additional order information at branch points.
For present purposes we do not need all this theory.

Proposition \ref{prop: Lorentz transversality} gives a sufficient
condition that can often be used to check
whether a section of a foliated circle bundle induces a slithering.

\section{Uniform foliations}
We will now specialize to the case of main interest: a slithering
$s:\Tilde M^m \to S^1$ of a
compact manifold $M^m$ around $S^1$.  Note that when 
when $\boundary M^m \ne \emptyset$, there is an induced
slithering of $\boundary M^m$. The foliation $\mathcal F(s)$ is
a codimension one foliation transverse to $\boundary M^m$. In particular
if $m = 3$ and $M^3$ is oriented, its boundary consists of tori.  

Any codimension one foliation $\mathcal F$ admits a transverse one-dimensional
foliation $\tau$ defined by any line field transverse to $\mathcal F$.
The pair of foliations gives a local $\R^{m-1} \times \R$
product structure for $M^m$.  For any parametrized arc
$\alpha:[0,t] \to M^m$ on a leaf of $\tau$ and any parametrized
path $p:[0,u] \to M^m$ on a leaf of $\mathcal F$, you can `comb' $\alpha$
along $p$ for some distance through the leaves of $\mathcal F$.
In other words,
there is a unique extension $\alpha \times p$ to a maximal
monotone subset $H$ of a rectangle,
satisfying
\begin{equation*}
\begin{align*}
\alpha\times p: H \subset [0,t] \times [0,u] &\to M  \\
\alpha \times p \restrict [0,t] \times 0 \;&=\; \alpha \\
\alpha \times p \restrict 0 \times [0,u]\; &=\; p \\
\bigl( r_1 \le r_2 \;\&\; s_1 \le s_2\bigr) &\implies \bigl ( (r_2,s_2) \in H \implies 
(r_1,s_1) \in H \bigr )
\end{align*}
\end{equation*}
where $\alpha \times p$ maps the two
coordinate directions to leaves of the two foliations.
In general, $H$ is an open set containing a neighborhood of the two
original sides, but not the full rectangle, because
the length of $\alpha$ might get longer
and longer, whipping out of control and failing to converge in the limit.
When the combing is interpreted as partially defining an
action of the groupoid of paths along leaves of $\mathcal F$ on arcs transverse
to $\mathcal F$, it is called the \df{holonomy} of $\mathcal F$.

\begin{definition}
The foliation $\mathcal F$ is \df{regulated} by $\tau$
if the holonomy of every arc $\alpha$ exists for all time along
any path. It is \df{uniformly regulated} by $\tau$ if the lengths of the
images of any arc $\alpha$ under the holonomy of $\mathcal F$
are bounded, with a bound that depends only on $\alpha$.

A foliation $\mathcal F$ is \df{uniform} if every 
closed transverse curve is non-trivial in homotopy, and if for every
pair of leaves $L_1$ and $L_2$ in the universal cover,
each is contained in a bounded neighborhood of the other.

Two uniform foliations $\mathcal F$ and $\mathcal G$ of a manifold $M$
are \df{uniformly equivalent} if for every pair of leaves $L$
of $\Tilde {\mathcal F}$ and $L'$ of $\Tilde {\mathcal G}$, each is
contained in a bounded neighborhood of the other.
\end{definition}
The prohibition on null-homotopic closed transversals in uniform
foliations eliminates examples that have a very different flavor,
such as the Reeb foliation (or any foliation) of $S^3$.
The condition implies that every leaf is properly
embedded in the universal cover, since a leaf in the universal cover
can never intersect a transverse arc more than once. In dimension $3$,
by the celebrated work of Novikov \cite{Novikov:foliations}, 
a transversely oriented foliation on any orientable manifold other than
$S^2 \times S^1$ satisfies this condition if and only if it does not
contain a Reeb component.

It follows from the definition that when
$\mathcal F$ is regulated by $\tau$, the lifts
of $\mathcal F$ and $\tau$
to the universal cover of $M^m$ define a product structure. Conversely,
if the leaves of $\mathcal F$ and $\tau$, lifted to the universal cover,
are the factors in a product structure, then $\tau$ regulates $\mathcal F$.
The product structure gives two slitherings for $M^m$---a slithering of $M^m$
around the universal cover of any leaf of $\mathcal F$,
and another slithering around $\R$ (which is the universal covering of
a leaf of $\tau$.) The two foliations complement each other,
serving as flat connections for the two slitherings.

When $\mathcal F$ is regulated by $\tau$, then 
it has no null-homotopic closed transversals; if the regulation
is uniform, then $\mathcal F$ is uniform.  If $M^m$ is compact and
$\mathcal F$ is uniform, it easily follows that any $\tau$ that
regulates $\mathcal F$ regulates it uniformly.  This is not true
for noncompact $M^m$ (\emph{e.g.} an easy counter-example can be
constructed on $\R \times S^1$.)

\medskip
In \cite{MR88m:58024}, \'Etienne Ghys gave an
elegant description of a certain equivalence relation on foliated circle bundles
in terms of bounded cohomology, and characterized equivalence classes in
terms of blowing-up of leaves.  This equivalence relation is a special case
of uniform equivalence of uniform foliations. Ghys' characterization
of equivalence classes, upon translation to the present context,
generalizes to the following:

\begin{proposition} \label{proposition: uniform equivalence}
Let $\tau$ be a 1-dimensional foliation of a compact $n$-manifold
$M$, and let $\mathcal F$ and $\mathcal G$ be uniformly equivalent
uniform foliations of $M$ that are regulated by $\tau$.

If every leaf of $\mathcal F$ and of $\mathcal G$ is dense,
then $\mathcal F$ is topologically equivalent to $\mathcal G$.

In any case, there is a third uniform foliation $\mathcal H$
that can be obtained from each of $\mathcal F$ and $\mathcal G$
by blowing up at most a countable set of leaves, where each blown-up leaf
is replaced by a foliated $I$-bundle.
\end{proposition}

\proof
Let $X$ and $Y$ be the spaces of leaves of $\Tilde {\mathcal F}$ and
$\Tilde{\mathcal G}$;
each is homeomorphic to $\R$. Let $I \subset X \times Y$ be the closure
of the set of pairs of leaves from the two foliations
that intersect each other. Each leaf of one foliation intersects
an interval's worth of leaves in the other, since all leaves
separate $\Tilde M$ into two components, so the intersection of $I$ with
any line $x \times Y$ or $X \times y$ is a compact interval (possibly
a single point.) Therefore, $\boundary I$ is the union
of two embedded lines (but the lines might intersect each other.)
Choose one of the lines, call it $l$. Define a submanifold
$M_l \subset \Tilde M \times \Tilde M$ as a union of
copies of $l$, one copy for each leaf of $\Tilde \tau$---this
makes sense because any leaf of $\Tilde \tau$ 
is canonically homeomorphic to $X$ and to $Y$. 

Observe that $M_l$ has the product structure of a leaf of
$\Tilde {\mathcal F} \times l$; it is invariant by the action
of $\pi_1(M)$, so the quotient $M_l / \pi_1(M)$ is homeomorphic
to $M$, and has a codimension one foliation $\mathcal H$. Projection
to the two factors shows that $\mathcal H$ is a blow-up of
$\mathcal F$ and $\mathcal G$

If leaves of $\mathcal F$ and of $\mathcal G$ are dense, then the
blowing up is trivial, since  one of the two projected images of
any blown-up region in $\mathcal H$ is a proper open invariant subset
for one of $\mathcal F$ or $\mathcal G$.
\endproof

Given a codimension one foliation, if the induced foliation of
its universal cover has a 
product structure, the foliation is called \df{ $\R$-covered}.
The following example shows that not
all $\R$-covered foliations have uniform spacing of leaves:

\begin{example} \label{example: Anosov regulation}
Let $\phi: T^2 \to T^2$ be an Anosov diffeomorphism of the torus,
let $T_\phi$ be its mapping torus, let $\mathcal F_s$ be the stable
foliation of $T_\phi$, and $\mathcal F_{uu}$ be the strong unstable
manifold. The universal cover of $T_\phi$ is $\R \; \times $
the universal cover of $T^2$. The foliations 
 $\mathcal F_{uu}$ and $\mathcal F_s$ can be represented in
$\R^3$ by a family of parallel lines and an orthogonal family of
parallel planes, so $\mathcal F_{uu}$ regulates $\mathcal F_s$.
However, $\mathcal F_s$ is not a uniform foliation.
\end{example}

However, this example seems to be fairly exceptional. It is hard to construct
$\R$-covered foliations on `generic' $3$-manifolds. This is partly
because it is hard to know the space of leaves in the universal cover,
but it seems likely that there is also a fundamental
obstruction.
\begin{conjecture} \label{conjecture: R-covered}
A foliation of a closed hyperbolic 3-manifold is $\R$-covered
if and only if it is a uniform foliation.
\end{conjecture}

Actually, it is remarkable that hyperbolic manifolds can
have any kind of $\R$-covered foliations at all, since
surfaces in hyperbolic space `want' to separate from each other and go off
in different directions;
most constructions of foliations on hyperbolic manifolds
yield foliations which are not $\R$-covered.
Perhaps it shouldn't be surprising if
$\R$-covered foliations on hyperbolic manifolds turn out, as conjectured,
to be quite special.

A rough rationale for this conjecture is that when a foliation
of a hyperbolic manifold is not uniform, there tends to be recursively
nested spreading of the leaves that
forces nearby leaves to limit to the sphere at infinity in topologically
separated ways.  This is related to
section \ref{section: pA}, which constructs a transverse pseudo-Anosov
flow that controls the geometry of leaves of a uniform foliation,
and also to section \ref{section: Peano curves},
which analyzes how leaves of uniform foliations limit to the sphere at infinity,
in continuous, sphere-filling curves.  These structures seem likely
to occur for any $\R$-covered foliation, and they give a fairly
explicit conjectural picture of what $\R$-covered foliations should look
like.  See section \ref{section: Anosov} for a discussion of one class of
$\R$-covered foliations that do turn out to be uniform.

\begin{proposition} \label{proposition: slithering uniform}
If $s$ is a slithering of $M^m$ around $S^1$, then
$\mathcal F(s)$ is a uniform foliation.
\end{proposition}
\proof

We can define a rough measure of separation between two leaves for a
slithering $s$, as follows.
Let $\Tilde s: \Tilde M^m \to \R$ be a lift of $s$ to a fibering
over $\R$, where $S^1 = \R / \Z$.
The fibers of $\Tilde s$ are connected.
Let $L_r$ and $L_t$ be any two fibers of $\Tilde s$, that is,
leaves in $\Tilde M$, where $r, t \in \R$. The interval
$[r,t]$ wraps some whole number of times around $S^1$, with some bit
left over. Define a function $z(L_r, L_t)$ to be the even number
$2 (t- r)$ when $t-r$ is an integer, and the odd number
$2\floor{t-r} + 1$ otherwise.
With this definition,
covering transformation of $\Tilde M \to M$ preserve the function $z$
on pairs of leaves. 
Similarly, we define $z(\alpha)$ for any path $\alpha$ in $M$
by taking any lift of $\alpha$ to $\Tilde M$,  and evaluating
$z$ on the leaves of its endpoints.

The \df{$z$-diameter} of $M$ is the maximum, over all pairs of
points $x,y \in M$, of 
the minimum value of $z(\alpha)$, where $z(\alpha ) > 0$,
$\alpha(0) = x$ and $\alpha(1) = y$.
Since $M$ is compact, its $z$-diameter is some finite number $k$.
Any arc $\alpha$ in $M$ such that $z(\alpha) > k$ must intersect every
leaf of $M$. In fact, if every leaf of $\mathcal F(s)$ is dense, then
given $x$ and $y$, the leaf through $y$ intersects any transverse
arc starting at $x$, so the $z$-diameter of $M$ is $1$.

Let $L$ and $L'$ be any two leaves in $\Tilde M$; assume that
$z(L,L') > 0$. Let $\beta$ be any arc transverse to $\mathcal F(s)$
with $z(\beta) = z(L,L') + k + 2 $. We can subdivide
$\beta = \beta_1 * \beta_2$, where $z(\beta_1) > k$ and
$z(\beta_2) > z(L, L')$.  

For any point $ \Tilde x \in L\subset \Tilde M$, we can
project to $M$, connect the image point $x$ to the $\beta_1$ by a path $p$
on its leaf, then lift $p$ and $\beta$
back to an arc $\Tilde \beta $ in $\Tilde M$ that intersects
$L$ in $\Tilde \beta_1$.  Since $z(\beta_2) > z(L,L')$, it follows that
$\Tilde \beta$ intersects both $L$ and $L'$.
The corresponding
lift of $\beta$ intersects $L'$, giving an upper bound
to the distance from $L$ to $L'$.  

By symmetry, there is also a bound when $z(L',L) < 0$---we can, for
instance, just reverse the orientation of $\R$.
\endproof
\begin{corollary} \label{corollary: rough distance}
The distance between any pair of leaves $L$ and $L'$ in $\Tilde M^m$
is bounded above by some constant times $|z(L, L')|$, and bounded
below by some constant times $|z(L, L')| - 1$.
\end{corollary}

There is a construction going in the reverse direction,
from uniform foliations to slitherings, but it is not
an exact converse.
We will give a statement expressed in terms of
a foliation that is uniformly regulated by a line field. The same proof
can be applied to uniform foliations in general, but in this
setting the conclusion would be weaker. This is
a variation of proposition \ref{proposition: uniform equivalence},
where one is constructing a
uniform equivalence from a foliation to itself:

\begin{theorem} \label{thm: uniform separation}
Let $\mathcal F$ be a codimension-one foliation of $M^m$ that is
uniformly regulated by a 1-dimensional foliation $\tau$.

If every leaf of $\mathcal F$
is dense then $M^m$ is the foliation of
a slithering of $M^m$ around $S^1$.

In any case, there is a slithering $s: \tilde M^m \to S^1$
of $M^m$, regulated by $\tau$ such that
$\mathcal F(s)$ is uniformly equivalent to $\mathcal F$.
\end{theorem}
\proof
This could be proven using the same technique as for proposition
\ref{proposition: uniform equivalence}, but we'll express the
proof in somewhat different language instead.
\begin{figure}[hbtp]
\centering
\begin{minipage}{.37\textwidth}
\includegraphics[width=\textwidth]{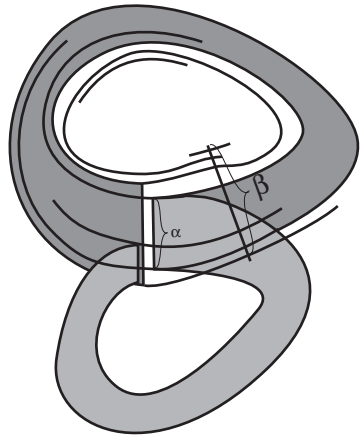}
\end{minipage}
\begin{minipage}{.63\textwidth}
\caption[Expanding calipers]
{\label{figure: calipers}
In a uniform foliation, we can start with any transverse arc $\alpha$,
and look at all its images under the holonomy of the foliation.
The limiting arc $\beta$ cannot expand any wider than it already
is. It acts like a pair of calipers that can roughly but consistently
measure the transverse  distances between leaves, since
holonomy images of $\beta$ can never nest with each other.}
\end{minipage}
\end{figure}

First we analyze the case that every leaf of $\mathcal F$ is dense.
Let $\alpha$ be any arc of $\tau$ (see figure \ref{figure: calipers}.)
Choose a Riemannian metric on
$M^m$ (or just a path-metric with reasonable properties,
if the data isn't very smooth).
Let $A$ be the supremum of the length of holonomy images
of $\alpha$, and let $\beta$ be any arc of length $A$ that is a
limit of images of $\alpha$. 

Then no holomy image of $\beta$ can have length greater than $A$, since
every holonomy image of $\beta$ is a limit of images of $\alpha$.
In particular, no holonomy image of $\beta$ can properly contain $\beta$.

Let $l$ be any flow line of $\tau$ in the universal cover, and
consider all the lifts of images of $\beta$ to $l$.  Since leaves
are dense in $M^n$, images of both endpoints of $\beta$ are dense.
The map that takes each lower endpoint of an image of $\beta$ to
its upper endpoint is a well-defined, monotone function 
from one dense subset of $l$ to another dense subset of $l$,
with monotone inverse.
Therefore, it extends continuously to a homeomorphism
$f:l \to l$.  This homeomorphism commutes with the action of $\pi_1(M^n)$.

In other words, the projection of $\Tilde M^n$ to $l$ modulo $f$
is a slithering around $S^1$.

\medskip
If not every leaf of $\mathcal F$ is dense, we can still carry out
most of the argument. Begin with an arc $\alpha$
of $\tau$ with both endpoints on a minimal set $K$ of $\mathcal F$.
We obtain a limiting
arc $\beta$ that has endpoints on $K$ and whose
holonomy images cannot nest with itself. If $l$ is a flow line in the
universal cover, then we obtain a monotone function $f$ with a monotone
inverse from $K \cap l$ to itself. Therefore $f$ is a homeomorphism
of $K \cap l$.

The only possibilities for a proper minimal set such as $K$ in a
codimension one foliation is that $K$ is
either a closed leaf, or an exceptional
minimal set (one where $K \cap l$ is a Cantor set.)

If $K$ is a closed leaf $L$, then $M$ actually fibers over $S^1$ with fiber
$L$ and structure group $f$; the foliation by fibers of the
fibration is uniformly equivalent to $\mathcal F$.

If $K$ is an exceptional minimal set, then we can collapse each arc
of $\tau \cap  M\setminus K $ to obtain a uniform foliation where
every leaf is dense, which therefore comes from a slithering.
\endproof

In high dimensions, one can modify a slithering by taking
the connected sum with a simply-connected manifold on each leaf.
Sometimes the result is a foliation whose leaves are not homeomorphic:
for instance, we could start with a 5-manifold with a slithering
similar to example \ref{example: torus slithering}, with a transverse
curve that does not intersect every leaf, then perform the leafwise
connected sum with $\CP^2$ along the curve.  This indicates the importance
of the condition that the foliation of the universal cover is a fibering.
It would appear that a variation of this procedure could yield a
manifold having a uniform foliation, but no slithering at all.

\begin{example}
Consider a foliated trivial $I$-bundle over a closed manifold with the
top glued to the bottom by a diffeomorphism $\phi$.
If there is a homeomorphism $h: I \to I$ that conjugates the holonomy
of the bundle to the holonomy composed with $\pi_1^*(\phi)$, then the
resulting foliation is the foliation of a slithering over $S^1$,
where the structure map $f$ for the slithering is constructed by
stringing together copies of $h$.  Otherwise, if the holonomy is not
invariant at least by some power of $\pi_1^*(\phi)$, the foliation
is not the foliation of a slithering.

This and other similar examples show that
foliations constructed by blowing up leaves of $\mathcal F(s)$
are not typically foliations of slitherings around $S^1$, although
they still slither around $\R$.
The structure map $Z$, which comes from a
generator of the group of deck transformation of $\R \to S^1$, is a
leaf-preserving homeomorphism isotopic to the identity in \(M\) but usually
not isotopic to the identity on a leaf.

Blowing-up operations for slitherings can be naturally performed in terms of
the foliated circle bundle over $M$,
rather than directly in terms of the foliation.
\end{example}

\subsection{Uniform regulation by Lorentz structures}

\label{subsection: Lorentz structures}
Let $\mathcal F$ be a foliation of a compact Riemannian
manifold $M^m$ that is uniformly regulated by a line field $\tau$.
From theorem \ref{thm: uniform separation} we see that there is
some constant $A$ such that given any two leaves $L$ and $L'$
in $\Tilde M$,
there is a sequence of intermediate leaves $L = L_0, L_1, \dots, L_n = L'$
such that the distance between $L_i$ and $L_{i+1}$ is uniformly
bounded by $A$. If $\tau'$ is another
line field that makes a sufficiently small angle with $\tau$, the flow
lines of $\tau$ and $\tau'$ stay reasonably close to each other by the
time they go a distance of $A$.  In particular, we can guarantee
that in the universal cover,
the flow lines of $\tau'$ hit $L_{i+1}$ in a distance only modestly
greater than $A$ after they hit $L_i$, or in other words,
$\mathcal F$ is also uniformly regulated by $\tau'$.

Estimates for this kind of information can often be conveniently encoded
by a Lorentz structure.
This can be done with a Lorentz metric, that is, a quadratic form $q$
of signature $(n,1)$, where $\tau$ is contained in the double cone
where $q$ is negative.  More generally, we could encode the information
with an open convex cone in the tangent space of each point (not necessarily
a quadratic cone). We'll call this a Lorentz cone structure. A Lorentz
cone structure is \df{transverse to $\mathcal F$} if every line contained
within the cone is transverse to $\mathcal F$. We say that a foliation
$\mathcal F$ is uniformly regulated by 
a transverse Lorentz cone structure $C$ if
any two leaves $L$ and $L'$ in $\Tilde M$ can be connected by a
transverse arc within $C$, and there is a finite upper bound $K(L,L')$
for the length of any arc within $C$ connecting
$L$ to $L'$.

As a limiting case, we will say
that a Lorentz cone structure $C$ \df{almost uniformly regulates}
$\mathcal F$ if
it is the increasing union of Lorentz
cone structures that uniformly regulate $\mathcal F$.  As an example,
consider the foliation by fibers of
any actual fibration $M^n \to S^1$, with the Lorentz
cone structure $C$ which is the union of the two open half-spaces that
exclude the tangent spaces to the fibers. Then $C$ almost uniformly
regulates the foliation.  In other examples just as in this one,
it is often easiest to describe and
think about a limiting case that almost uniformly
regulates $\mathcal F$.

If $M$ is a manifold with a Lorentz cone structure $C$ and
$g: N \to M$ is a differentiable map, we'll say that $g$ is
\df{transverse to $C$} if for each $x \in N$ there is a tangent
vector $V \in T_x(N)$ taken into $C$, $dg(V) \in C$.
In that case, for each open
convex half $H$ of the double-cone $C$, the set of vectors
that map to $H$ form a convex cone in $T_x(N)$,  describing
a Lorentz cone structure $g^*(C)$.  We can think of a foliation
as a special case of Lorentz cone structure, where each open convex
half-cone is a half-space; this definition is a generalization
of the definition
of a map transverse to a foliation, and of the pullback foliation
$g^*(\mathcal F)$.
\begin{proposition} \label{prop: Lorentz transversality}
If $\mathcal F$ is a codimension one foliation of $M$ almost uniformly
regulated by a Lorentz cone structure $C$, and if $g: N \to M$
is a differentiable map transverse to $C$, then $g^* \mathcal F$ is
a codimension one foliation almost uniformly regulated by $g^* C$.
\end{proposition}
\proof
This follows from compactness considerations: if $C' \subset g^*C$
is a Lorentz cone structure whose half-cones have closure contained
in $g^* C$, then the images $g(C')$ have closure contained in $C$.
\endproof

It is worth observing that this proposition can be used to give a compact
geometric criterion for slitherings in terms of their
associated foliated circle bundles.  Any foliated circle bundle
can be almost uniformly regulated by a Lorentz cone structure $C$,
and in particular cases there are explicit constructions using for
example curvature estimates. Given a foliated circle bundle $E \to M$,
a section $S \to E$ induces a slithering if it is transverse to some
$C$.  

\begin{figure}
\centering
\includegraphics[width=\textwidth]{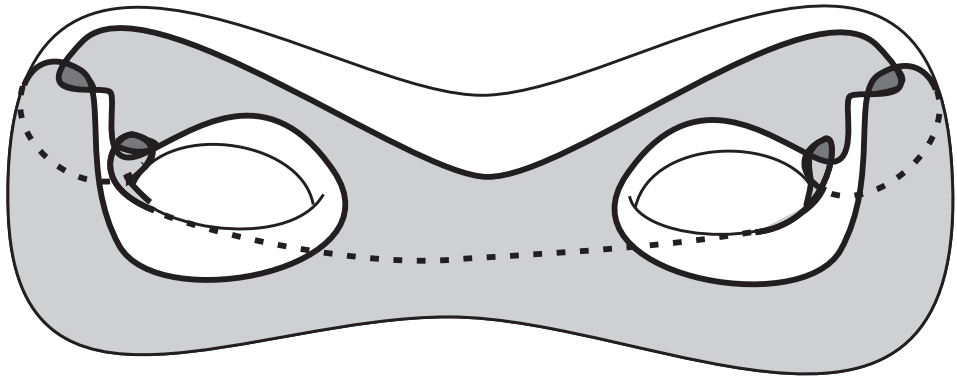}
\caption[A timelike knotted curve in $T_1(M^2)$]
{\label{figure: timelike knot}
The tangent to this curve is a time-like curve for the
Lorentz structure of $T_1(M^2)$ for certain hyperbolic structures on $M^2$.
The curve is homotopic to a fiber, since it is the boundary of an immersed
disk. An annulus realizing the homotopy in $T_1(M^2)$ intersects itself in
short double arcs contained in the regions of self-intersection on the surface.
The image of the disk generates $\pi_1(M^2)$, from which one can deduce that
the complement of the curve is atoroidal.
There are hyperbolic metrics that make the disk round,
so in the universal cover view of figure \ref{figure: circle at infinity}
the curve is covered by a countable set of helices,
each interlinked with four neighbors.
The complement of this curve (or any other time-like
curve) slithers around $S^1$, as does every cover of $T_1(M^2)$
branched over timelike curves.
}
\end{figure}
\begin{example}

Let $M^2$ be a hyperbolic surface, and $\mathcal F$ the circle at
infinity foliation of $\TS(M^2)$ (example \ref{example: tangent of hyperbolic}.)
The tangent to a horocycle in $\Hy^2$
almost traverses the circle at infinity, omitting
only one point, where the horocycle is tangent to $S_\infty^1$.  Similarly,
a curve in $\TS(M^2)$ whose base point follows a horocycle in $\Hy^2$, lifted
to a vector that makes a constant angle to the tangent to the
horocycle, but doesn't point to the point of tangency on $S_\infty^1$,
traverses the entire circle except for that one point.  All such curves
together sweep out the boundary of a Lorentz cone structure for $\TS(M^2)$,
where a curves in $\TS(M^2)$ are inside the cone if 
their vectors turn faster than their base points move.
One can think of time-like trajectories as dancers
moving in $\Hy^2$ in a way that all the scenery, near and far, in front
and behind, appears
to rotate consistently in one direction.  This Lorentz cone structure
comes from a Lorentz metric given by the Killing form on $\PSL(2,\R)$.
This Lie group has the same complexification ($\PSL(2,\C)$) as $SO(3)$,
and its Lorentz metric lifted to $\SL(2,\R)$ has analytic continuation
that agrees with the round metric on $S^3$.

\begin{itemize}
\item
If $Q^2$ is any closed hyperbolic orbifold and
we remove any closed time-like trajectory from $\TS(Q^2)$,
the resulting manifold still slithers around $S^1$.  There are
many possible closed time-like trajectories.
For example, the tangent to any curve with curvature greater than 1
is a time-like trajectory; it gives an embedded curve in $\TS(Q^2)$
as long as it is never tangent to itself.

Every homotopy class of curves in $Q^2$
has infinitely many regular homotopy classes of immersions
that can be arranged to satisfy this condition. In fact, all but one of
the $\Z$'s worth of regular
homotopy classes can be made time-like.  This can be done by
transporting a large-diameter circle immersed in $Q^2$ so that its center
traverses a geodesic in the base. A pencil speeds around the circle
drawing an immersed curve as the circle moves. The pencil
travels at constant velocity in parallel-translated coordinates for the
circle.

Furthermore, each time-like regular homotopy class is represented by
many distinct time-like knots.  In the case of
tangents of immersed curves of curvature $> \; 1$,
the knot type usually changes whenever a tangency between the curve and itself
occurs. For the moving circle construction, infinitely many events
of this type occur as the circle's radius tends to infinity.
\item
If $M^3$ slithers around $S^1$ and $\mathcal F(s)$ is almost uniformly
regulated by a Lorentz cone structure $C$, then
any branched cover of $M^3$ over any time-like link
also slithers around $S^1$.  
\item
We can remove any time-like curve from $M^3$ as above,
and replace it by any 3-manifold
with torus boundary that fibers over $S^1$, to obtain a new
manifold that slithers around $S^1$.  
\item
Let $\tau \subset M^3$ be a $1$-dimensional train-track embedded transversely
to $C$. Suppose that we have assigned integral weights to the branches
of $\tau$ in a way that satisfies the switch additivity condition. Then
we can remove a neighborhood of $\tau$, and for any branch labeled $g$,
insert a surface of genus $g$  $\times [0,1]$.  We have various choices
of gluing maps to glue the ends of the units together at each switch.
These can be arranged, if desired, so that the entire inserted
assemblage has the same homology as the neighborhood of $\tau$ that was
removed.  The resulting manifold slithers around $S^1$, since it has a
map to $M^3$ that is transverse to $C$.
\end{itemize}
\end{example}

In light of the fact that some Seifert fiber spaces are homology spheres,
and the fact that there are pseudo-Anosov maps of surfaces that
induce every possible symplectic automorphism of first homology, 
these constructions are powerful enough to produce an atoroidal
3-manifold that slithers around $S^1$ with any desired 
homology type. However, these constructions only scratch the
surface, since in these examples
the action of $\pi_1(M^3)$ on $S^1$ still factors through a Seifert fiber
space. The moduli of representations $\pi_1(M^3) \to \widetilde{\Homeo (S^1)}$
seems to be quite interesting---even when limited to representations
in restricted subgroups, such as
$\widetilde {\PSL(2,\R)}$ or
or PL-homeomorphisms---and deserves investigation
beyond the present scope.

\section{Groups, inequalities and topology}
\label{section: groups, inequalities and topology}

The Milnor-Wood inequality for foliated circle bundles and Stalling's
characterization of fundamental groups of 3-manifolds that fiber over
$S^1$ demonstrate an interplay of group theory, geometry and
topology that is involved in group actions
on $\R$ and on $S^1$, and in manifolds that slither around $\R$ or
$S^1$. This section will recount some of the rudimentary theory of
this interaction. 

One theme is that groups of periodic homeomorphisms have
approximate homomorphisms to $\R$, despite the fact that they might not
have any actual homomorphisms. This theme can be expressed
in several ways.

\begin{figure}
\begin{minipage}{.35\textwidth}
\includegraphics[width=\textwidth]{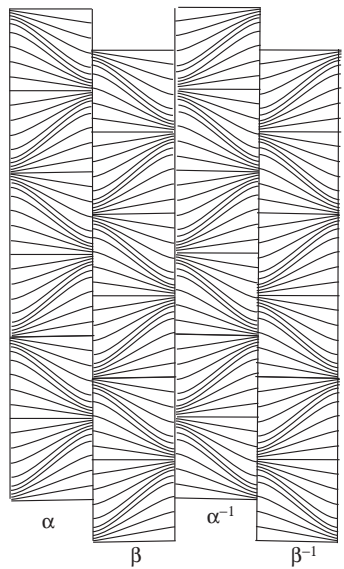}
\end{minipage}
\begin{minipage}{.65\textwidth}
\caption[Commutators can translate]
{ \label{figure: commutators can translate}
This diagram shows the commutator of two periodic diffeomorphisms $\alpha$
and $\beta$ of $\R$, where $\alpha$ and $\beta$ both have fixed points,
but their commutator translates every point downward, sometimes
further than one period. 
By adjusting $\alpha$ and $\beta$, in the limit some points go downward
nearly two periods, showing that proposition
\ref{proposition: commutator bound} is sharp.
}
\end{minipage}
\end{figure}

\begin{proposition}\label{proposition: commutator bound}
Let $a, b, c \in \Homeo_+(\R)$ satisfy $[a,c] = [b,c] = 1$.
Suppose that for all $x \in \R$, $c(x) > x$. Then, for all
$x \in \R$, $c^{-2} (x) < [a,b](x) < c^2(x)$.
\end{proposition}
See figure \ref{figure: commutators can translate} for a picture
of what this is all about.

\proof
For any $x\in \R$ there are integers $k$ and $l$ such that
\begin{align*}
c^k(x) \le a(x) &< c^{k+1}(x) \\
c^l(x) \le b(x) &< c^{l+1}(x)
\end{align*}
Therefore 
\begin{gather*}
c^{k+l}(x) \le a( c^l(x)) \le a(b(x))  < a(c^{l+1}(x)) < c^{k+l+2}(x) \\
c^{k+l}(x)  \le  b(c^k(x))  \le b(a(x)) < b(c^{k+1}(x)) < c^{k+l+2}(x)  \\
\end{gather*}
Let $x' = b(a(x))$, so $x = a^{-1}b^{-1} x'$. Express the second
line of inequalities
in terms of $x'$, then apply the first, to obtain:
\begin{gather*}
c^{-k-l-2} (x') < a^{-1}b^{-1}(x')  \le c^{-k-l} ( x') \\
 c^{-2}(x') <  a b a^{-1} b^{-1} (x') < c^2(x')
\end{gather*}
Since $x'$ is arbitrary, this proves the proposition.
\endproof

\begin{proposition} \label{proposition: groups w/o fixed points}
If $G \subset \Homeo(\R)$ acts without fixed points (that is, no
element of $G$ except the identity fixes any point in $\R$),
then $G$ is isomorphic as a group to a subgroup $R(G)$
of the additive group of $\R$, with action obtained from
$R(G)$ by blowing up at most a countable set of orbits.
\end{proposition}
(See also \ref{proposition: space-like generators}.)
\proof
If $G$ acts without fixed points, then $G$ has a linear ordering
where for $a,b \in G$, $a < b$ if for all $x \in \R$, $a(x) < b(x)$.
This linear ordering is bi-invariant, that is, if $a < b$
then for all $g \in G$, $ga < gb$ and $ag < bg$.

If there is a least element of $G$ that is greater than $1$, then
$G$ is cyclic, and we are done.

Otherwise, given any element
$c \in G$ (we imagine $c$ to be small), we can compare the
commutator of $a $ and $b$ with powers of $c$; the same inequalities as in
proposition \ref{proposition: commutator bound} and its proof hold,
but justified now by the bi-invariance of order rather than
commutativity.  In other words, the commutator of any two elements
cannot be greater than the square of any positive element. But there
is no lower bound to such squares, because
for any $1 < c \in G$, if $1 < c_1 < c < c_1^2$,
then $(c c_1^{-1})^2 <  (c c_1^{-1}) (c_1) = c$.  Therefore $[a,b] = 1$.

The rest of the statement follows by straightforward reasoning.
One method is to construct a measure on $\R$ invariant by $G$.
The integral of the measure gives a semi-conjugacy of the
action to an additive subgroup of $\R$---in other words, the
action of $G$ is obtained by blowing up at most countably
many orbits.
\endproof

We can apply this information to describe the centralizer of
a group of periodic homeomorphisms of $\R$.
We'll use coordinates
$S^1 = \R/\Z$, so for all $x \in \R$, $g(Z(x)) = Z(g(x)$ where $Z(x) = x+1$.

The multiplicative semigroup of $\Z$ acts by conjugation on such
representations. That is, if $n \ne 0$, we can define a new action 
$n \times \rho (x) = (1/n) \rho( n x)$.  This gives a new action
centralized by the supergroup of finite index $(1/n) \Z$.

\begin{proposition} \label{proposition: periodic centralizer}
Let $G$ be any group and $\rho: G \to \widetilde{\Homeo_+(S^1)}$ a
homomorphism such that all orbits for the action of $\rho(G)$ on $\R$
are dense.
Then $\rho$ is $\Z$-equivariantly conjugate to an action whose
centralizer is a closed subgroup of the additive group of $\R$.
Either $\rho(G)$ is abelian and isomorphic to an additive subgroup of $\R$,
or $\rho$ is is conjugate to $n \times \rho'$ where
$n \ge 1$ and the centralizer of $\rho'(G)$ is $\Z$.
\end{proposition}
\proof
Let $Z(\rho(G)) \subset \Homeo(\R)$ denote the 
centralizer of $\rho(G)$. For any element $h \in Z(\rho(G)$,
the fixed point set of $h$ is invariant by $G$; since all orbits of
$G$ are dense, this implies that either $h$ is the identity, or it has no
fixed points.  Therefore, $Z(\rho(G))$ is a subgroup of $\R$ containing $\Z$.
If it is cyclic, then $\rho$ is conjugate to a group
of the form $n \times \rho'$. (We find this conjugacy by looking at
the quotient circle $\R / Z(\rho(G))$.) Otherwise,
$Z(\rho(G))$ is isomorphic to a dense subgroup of $\R$. This group
cannot have an exceptional minimal set, since $G$ has dense orbits,
so the action of $Z(\rho(G))$ on $\R$ is standard.
\endproof

Similarly, if $s: \Tilde M \to S^1$ is a slithering, we can define
$n \times s$ by composing with the covering map $S^1 \to S^1$ of
degree $n$.  We'll say that a slithering is \df{primitive} if it is
not isomorphic to $n \times s'$ for any $n > 1$.

\begin{corollary} \label{corollary: slithering center}
Let $s: \Tilde M^3 \to S^1$ be a primitive slithering of a compact $3$-manifold
such that every leaf of $\mathcal F(s)$ is dense. 
For simplicity, assume also that $M^3$ is orientable and $\mathcal F(s)$
transversely orientable.
Let $\rho: \pi_1(M^3) \to \widetilde{\Homeo(S^1)}$ be the
action of its fundamental group. 

If $\pi_1(M)$ contains a non-trivial central element $a$,
then either
\begin{enumerate}
\item[i.] $M^3$ fibers over a torus, and $s$ is induced 
from a slithering of the base, or
\item[ii.] $s$ admits an invariant
measure and is a perturbation of a fibering
of $M^3$ over $S^1$, or
\item[iii.]
the action of $a$ satisfies
$\rho(a)(x) = x+l$ for some integer $l$.
\end{enumerate}
\end{corollary}
\proof
If $a$ is an element of the center of $\pi_1(M^3)$, then proposition
\ref{proposition: periodic centralizer} says that $a$ acts
either as the identity or without fixed points.

If $a$ acts as the identity and if
$\mathcal F(s)$ has no holonomy, then
$\rho(\pi_1(M^3)$ is abelian and admits an invariant measure, by
proposition \ref{proposition: groups w/o fixed points}, so it
fits in alternative (ii).

If $a$ acts as the identity and if $b$ is a loop on 
any leaf $L$ that has non-trivial holonomy,
then $a$ and $b$ generate a rank two abelian subgroup of $\pi_1(L)$;
$L$ is finitely covered by a torus, therefore it is a torus,
and $M^3$ fibers over $S^1$ with fiber $L$. The element $a$ is
invariant under the gluing map for this bundle, and so $M^3$
can also be described as fibering over $T^2$ with $a$ as fiber,
where the slithering respects this structure.

If $a$ does not act as the identity,
then \ref{proposition: periodic centralizer} gives us alternatives
(ii) or (iii).
\endproof

\subsection{Rotation numbers and commutator length}

There are a number of variations and extensions of the inequality
\ref{proposition: commutator bound}.  There is a rich literature
on this subject, which started with Milnor's (\cite{Milnor:inequality})
analysis of flat $\SL(2,\R)$-bundles.
Wood's paper (\cite{Wood:inequality}) analyzed circle bundles over surfaces by
studying $\widetilde{\Homeo_+(S^1)}$ and is
close to the present point of view.
Other noteworthy references include Eisenbud, Hirsch and Neumann
(\cite{EHN:SeifertFoliations}, analyzing Seifert fiber spaces), Sullivan
(\cite{Sullivan:EulerClass}, extending the theory to a more general context
and more general bundles) and Ghys (\cite{MR88m:58024}).  We will 
discuss this topic further in \cite{Thurston:circlesII} so for now
we will make just a quick foray, leaving most details, discussion  and
development to other sources.

When we have a homomorphism $\rho: G \to \widetilde {\Homeo(S^1)}$,
then for each $a \in G$ there is a \df{rotation number} $r(a) \in \R$,
where $r(a)$ generates the subgroup of $\R$ that best approximates
the cyclic subgroup generated by $a$. More precisely, $r(a)$ can
be characterized as the unique real number such that for all $x \in \R$,
the difference $|a^k(0) - k r(a)|$ has a bound independent of $k$.
The value $r(a) \mod \Z$ is the same as usual
rotation number of the action of $a$ as a homeomorphism of $S^1$; the
lifting to $\R$ is in virtue of the fact that a section is defined
for the associated torus bundle over $S^1$.   One way
to define and compute the rotation number is with
an invariant measure.  There is always
at least one measure $\mu$ on $\R$
invariant by $a$ such that $[0,1)$ has measure $1$.
Then we can define
$r(a) = \mu( [0, a(0)) ) = \mu([x, a(x)))$; the characterization
 of $r(a)$ by boundedness of $|a^k(0) - k r(a)|$ is
immediate.

An element $a \in G$ is \df{space-like} if $a$ has
a fixed point, or equivalently, if $r(a) = 0$. It
is a \df{positive time-like} element or simply \df{positive}
if $r(a) > 0$, or equivalently, $\forall x \; a(x) > 0$. 
Negativity is defined similarly. 

Using rotation numbers, we can amplify proposition
\ref{proposition: groups w/o fixed points},
for groups of periodic homeomorphisms of $\R$:
\begin{proposition} \label{proposition: space-like generators}
Let $G \subset \widetilde{\Homeo_+(S^1)}$ be any subgroup.
Then either
\begin{enumerate}
\item[(a)]
$G$ is generated by its space-like elements, or
\item[(b)]
the subgroup $G_0$ generated by space-like elements has
a common fixed point, and consists entirely of space-like elements.
This is equivalent to the existence of
an invariant measure on $\R$ for the action of $G$, and also
equivalent to the condition that $r$ (rotation number) is a homomorphism.
\end{enumerate}
\end{proposition}
\proof
The graph of a periodic homeomorphism of $\R$ can be mapped
to the cylinder $\left (\R \times \R\right ) / \Z$, where 
$\Z$ acts diagonally. The graph becomes
a closed curve that represents the generator of the fundamental
group of the cylinder. Thus,
periodic homeomorphisms are in 1--1 correspondence
with closed curves on the cylinder that are transverse to the
images of lines parallel to the axes. In rotated coordinates
$S^1 \times \R$, they become graphs of functions $S^1 \to \R$ that
strictly decrease the metric.

The set $S$ of space-like elements is represented by graphs that intersect
the diagonal $\Delta$, which is the graph of the identity. The
set $S*S$ of products of two space-like elements  is represented
by graphs that intersect these; $G_0$ is represented by graphs that
can be connected to $\Delta$ by a finite sequence of these curves,
such that each curve intersects the next in the sequence.

If you can go an unbounded distance 
in the cylinder by stepping from one curve  to intersecting curve, then
then the union of the graphs of elements of $G_0$ intersects every
possible curve in the homotopy class of the generator, so $G = G_0$,

Otherwise, let $A$ be the subset of the cylinder consisting of
the union of the graphs of elements of
$G_0$, along with all
bounded components of the complement.  Note that
$A$ is the graph of a symmetric,
transitive periodic relation on $\R$ (\emph{i.e.,} an equivalence
relation) modulo $\Z$.
From the defining properties, $A$ is invariant by $G_0$
and disjoint from its images by non-trivial elements of $G/G_0$.
These images are linearly ordered. If the linear order is discrete,
then $G/G_0$ is infinite cyclic. Otherwise,  the order-completion
is homeomorphic to $\R$---this $\R$ is the set of equivalence classes
of $\Tilde{\Bar A}$. Since
$G/G_0$ acts on $\R$ without fixed
points, we can apply \ref{proposition: groups w/o fixed points}.

For any $x_0 \in \R$, the intersection of line $y=x_0$ with $A$ is
the minimal interval containing its $G_0$-orbit;
 its upper and lower endpoints are necessarily fixed points
for the action of $G_0$.

\endproof

\begin{proposition} \label{proposition: uniform rotation}
Let $s_1, s_2: \Tilde M \to S^1$ be two slitherings of 
a compact manifold $M$ around
$S^1$. Then $\mathcal F(s_1) $ and $\mathcal F(s_2)$ are
uniformly equivalent if and only if the
rotation number functions for $s_1$ and $s_2$
agree up to a constant multiple.

If $\mathcal F(s_1)$ does not admit a
transverse invariant measure, then it is uniformly equivalent to
$\mathcal F(s_2)$ if and only if $s_1$ and $s_2$ determine the
same sets of space-like elements of $\pi_1(M)$.
\end{proposition}
\proof
One direction of the proposition is easy: if the foliations
are uniformly equivalent, then the
rotation number functions agree up to a constant multiple, since
the rotation number of $\alpha \in \pi_1(M)$ is defined by
the asymptotic translation distance of $\alpha^k$. The rotation
number function up to a constant
of course determines the set of space-like elements, which is the
set of elements whose rotation number is 0.

In the other direction, suppose that the two rotation number functions
$r_1$ and $r_2$
agree up to a constant multiple, $r_1 = C r_2$  Consider the map
$s_1 \times s_2: \Tilde M \to \R \times \R$. We claim that the image
of $s_1 \times s_2$ is contained in a bounded neighborhood of a line
in the plane. To establish this,
choose a compact fundamental domain $K$ for the action of $\pi_1(M)$
by deck transformations on $\Tilde M$. Let $z_i$ be the function on
pairs of leaves of $\Tilde {\mathcal F}(s_i)$ that roughly measures
twice the integer part of their separation.
Then for any $\alpha \in \pi_1(M)$ and any leaf $L_i$ of
$\Tilde {\mathcal F}(s_i)$ we have 
\[ 2 r_i(\alpha) - 1 \le z( L_i, \alpha(L_i)) \le 2 r_i(\alpha) + 1 .\]
Applying this to pairs of leaves that intersect $K$, we see that
the images $s_1 \times s_2 ( \alpha(K))$ are at a bounded distance from 
the line $y = C x$. By applying corollary \ref{corollary: rough distance},
we see that any
such pair of leaves has a uniformly bounded distance, hence all
pairs of leaves from the two foliations are contained in bounded
neighborhoods of each other.

If $\mathcal F(s_1)$ does not admit a transverse invariant measure,
then $\pi_1(M)$ is generated by space-like elements.
For any $\alpha \in \pi_1(M)$, define $g_i(\alpha)$ to be the minimum
length of its expression as a product of space-like generators,
and define
\[
h_i(\alpha)= \lim_n \to \infty \left ( g_i(alpha^n) \right )^{1/n}.
\]
By considering the cylinder discussed in the proof of proposition
\ref{proposition: space-like generators}, it is clear that $g_i(\alpha)$ can
be approximated by a constant multiple of $r_i$, up to a bounded additive
error term, so $h_i$ is a constant times rotation number. Therefore,
the set of space-like elements determines the uniform equivalence
class of $\mathcal F(s_i)$.
\endproof

\medskip
Inequalities such as in proposition \ref{proposition: commutator bound}
and its proof can be reworked and extended in terms of rotation numbers. The
rotation number of elements of a group $G$ of periodic homeomorphisms
can be thought of as a $1$-cochain on the group.%
\footnote{Cochains \emph{etc.}
for a group $G$ are the same thing as simplicial cochains \emph{etc.}
for the standard model of the Eilenberg-MacLane space
$K(G,1)$, whose $n$-simplices are labelings of the oriented $1$-skeleton of
an $n$-simplex by elements of
$G$ so as to form a commutative diagram.  A labeling is determined by its
value on a spanning tree, and different notations arise from different
choices of spanning tree.  We'll use the linear spanning tree
$\simplex{0,1}, \simplex{1,2}, \dots, \simplex{n-1,n}$. }
Proposition
\ref{proposition: space-like generators} describes the circumstances
that $r$ is a cocycle, which is equivalent to being a homomorphism.
Rotation number is usually not a cocycle, but its coboundary
$(\delta r)(a,b) = r(a) - r(ab) + r(b)$ is bounded:

\begin{proposition} \textbf{Milnor--Wood inequality for surfaces with boundary}
\label{proposition: Milnor-Wood with boundary}
For all $a, b \in G$,
\begin{gather}
\left |(\delta r)(a,b) \right | =
\left |r(a) - r(a b)+ r(b)  \right | \le 1 \label{eq: rot delta}
 \\
| r([a,b]) | \le 1 \label{eq: rot com 1} .
\end{gather}
Furthermore, if $a_1, b_1, a_2, b_2, \dots, a_n, b_n$
is any sequence of $2n$ elements
of $G$,
\begin{equation}
| r( [a_1, b_1] [a_2, b_2] \dots [a_n, b_n])| < n+1 , \label{eq: rot com 2}
\end{equation}
In general, for any homomomorphism of the fundamental group
of an oriented surface $M^2$ with boundary into $\widetilde{\Homeo_+(S^1)}$,
the sum of the rotation numbers of its boundary curves
does not exceed $\max(0,-\chi(M^2))$.
\end{proposition}
\proof
The various inequalities all
follow from the statement about an oriented surface with
boundary: (\ref{eq: rot delta}) is the case of pair of pants,
(\ref{eq: rot com 1}) is a punctured torus, and (\ref{eq: rot com 2})
is a punctured surface of genus $n$.

Here is a sketch of a proof. Given the representation of $\pi_1(M^2)$,
let $\xi\to M^2$ be the associated foliated circle bundle, with section
$S: M^2 \to \xi$ that is
defined up to homotopy by the structure group $\widetilde {\Homeo_+(S^1)}$.
Choose a complete hyperbolic metric of finite area on $M^2$, and lift this
to give a hyperbolic metric for the leaves of the foliation.
By \cite{MR84j:58099}, there is a harmonic measure for this foliation,
which can be normalized so it gives measure $2\pi$ to each fiber, which
assigns an $\SO(2)$ structure to each fiber, and makes the bundle
into a principle $\SO(2)$-bundle. The foliation gives a plane field
transverse to the fibers; average the slopes of the images of this
plane field under the action of the circle $\SO(2)$.  The resulting
connection has curvature that never exceeds $1$. Near the cusps,
it converges to a flat connection whose slope is $2 \pi$ times
the rotation number for the corresponding boundary component.

This implies that the integral of the curvature does not exceed
$- 2 \pi  \chi(M^2)$, which translates into the inequality that
the sum of rotation numbers of the boundary components is not greater
than the Euler characteristic of $M^2$.
\endproof
\begin{remark}
There are better inequalities in the cases when the
the rotation numbers are not all integers. For instance, if
$0 < r(a)  < 1$, then for all $x$, $a^{-1}(x) > x$, and 
the same holds for the conjugate: $b a^{-1} b(x) < x$. It follows 
that
\[
r([a,b]) = r( a * (b a^{-1} b^{-1})) < r(a),
\]
and similarly $r([a,b]) > r(a^{-1})$. If $1-r(a) < r(a)$, we can
get a better inequality by using the fact that
$a b a^{-1} b^{-1} = (a Z^{-1}) b (Z a^{-1}) b^{-1}$
and applying the same reasoning. Continuing along the same path,
as long as at least one of $r(a)$ or $r(b)$ is not an integer, we
can approximate it by the nearest power of $Z$, and deduce that
$r([a,b]) \le 1/2$. It is curious that when $r(a)$ and $r(b)$ are
both $0$ (or other integers), then $r([a,b])$ can attain
$\pm 1$, as shown by the the tangent bundle of a once-punctured
surface, or by the example of figure \ref{figure: commutators can translate}.

This is related to the fact that most periodic homeomorphisms with a
rational rotation
number have neighborhoods in the group of homeomorphisms where
the rotation number is constant.  In particular, the cases of integral
rotation numbers are big baskets.  In
$\widetilde{\PSL(2,\R)}$, these
baskets contain one lift for every hyperbolic element. 
\end{remark}

For an element $g$ of the commutator subgroup $G' = [G,G]$, the
\df{commutator length} $\cl(g)$ is the minimum number of commutators
$[a_i, b_i]$ whose product equals $g$. This is a sub-additive function:
$\cl(gh) \le \cl(g) + \cl(h)$.
When $g$ is not in $G'$, then we can define $\cl(g) = \infty$.
Define $\acl(g) = \liminf_{n\to\infty} 1/n \cl(g^n)$; we'll call
this the \df{asymptotic commutator length} of $g$; it is finite
if and only if $\alpha$ maps to an element of finite order in $H_1(M)$

\begin{corollary}
If $M$ slithers around $S^1$, then for any $\alpha \in \pi_1(M)$,
\[
r(\alpha) \le \acl(\alpha).
\]

A 3-manifold whose fundamental group contains a central
element $\alpha$ cannot slither around
$S^1$ unless $M$ is a nilmanifold or $\acl(\alpha) \ge 1$.

A 3-manifold whose fundamental group contains a normal $\Z^2$ subgroup
cannot slither around $S^1$ unless $M$ has a Euclidean or nilgeometry
structure.
\end{corollary}
This contains a version of the Milnor-Wood inequality expressed in terms
of slitherings: for a circle bundle of Euler class $n$ over $M^2$, the fiber
has order $n$ in homology. The asymptotic commutator length of
the fiber is $\chi(M^2)/n$, which must be at least $1$ to admit
a foliation transverse to the fibers.
\proof
The inequality on rotation numbers
follows from equation \ref{eq: rot com 2}. The application to 3-manifolds
whose fundamental group has a center follows from
\ref{corollary: slithering center}.

Rotation number is a linear function on any abelian subgroup of
$\pi_1(M)$: this follows directly from the definition, and also
is a consequence of inequality \ref{eq: rot delta} of
proposition \ref{proposition: Milnor-Wood with boundary}. If there is
a normal $Z^2$ subgroup $A \subset \pi_1(M^3)$, then the linear
function $\rho \restrict A$ must be invariant by conjugacy, if
$\mathcal F(s)$ is transversely orientable, otherwise it is invariant
up to sign. If
$M^3$ is orientable, this implies there is a common
eigenvector with eigenvalue identically $1$ for the conjugacy
action of $\pi_1(M^3)$ on $A$, \emph{i.e.},  there is a central
$\Z\subset \pi_1(M^3)$.

In the non orientable or non-transversely-orientable
cases, there still is a normal $\Z$ whose centralizer necessarily has
index at most two. It still follows that $M^3$ fibers with
fiber a circle with Euclidean base,  but besides the torus, the
base might be a Klein bottle ($(\mathsf{XX})$ in Conway's orbifold notation), or any
of the Euclidean 2-orbifolds whose non-trivial local groups have order $2$:
an annulus ($**$), a Moebius band ($*\mathsf X$), a pillow ($2222$),
a sack ($22*$) or a cross-sack ($22\mathsf X$).
\endproof

\section{Geodesic currents}
\label{section: geodesic currents}

The uniformization theorem says that
every Riemannian metric on a surface is conformally equivalent to
a metric of constant curvature. Alberto Candel (\cite{Candel:Uniformization})
analyzed how the uniformization varies from leaf to leaf in a lamination
with 2-dimensional leaves.  In particular, he showed that if $\mathcal F$
is a codimension one foliation of an irreducible 3-manifold
that either the foliation is a perturbation of a foliation with torus leaves,
or there is a Riemannian metric that has constant negative curvature on
each leaf.  This Riemannian metric varies continuously, but
may not be very differentiable in the manifold as a whole.  We will
not concern ourselves here with questions of regularity of the metric, since
we will really only use the quasi-isometric properties, which are not
delicate at all. 
If preferred, the metric of constant curvature $-1$
on each leaf can be approximated by a $C^\infty$ or even $C^\omega$
Riemannian metric on $M$ that induces
a metric on each leaf whose curvature is pinched
between $-1-\epsilon$ and $-1+\epsilon$;
such a metric would be more than adequate for what we will do.

Given any Riemannian metric for the leaves of a codimension one
foliation $\mathcal F$ of a closed $3$-manifold $M^3$, let $T_1(\mathcal F)$
denote the unit tangent bundle of the leaves of $\mathcal F$, and let
$\gfl: \R \times T_1(\mathcal F) \to T_1(\mathcal F)$ denote the geodesic
flow on the leaves.  

Every flow on a compact space admits at least
one invariant measure; usually, there are many different invariant
measures.  Here is a more explicit method to find 
invariant measures, in the case of $\gfl$ for the foliation
$\mathcal F(s)$ of a slithering $s$ of a compact $3$-manifold $M$.
Proposition
\ref{proposition: space-like generators} says that the normal closure of
the fundamental groups of the leaves of $\mathcal F(s)$ is the kernel
of a homomorphism (typically the trivial homomorphism)
of $\pi_1(M)$ to a subgroup of $\R$. The image
group consists of periods (integrals around closed curves)
for a transverse invariant measure equipped with a transverse orientation,
turning it into a closed current (which is a generalization of a
closed $1$-form.)  If $\pi_1(M)$ is not abelian, the kernel is automatically
non-trivial, therefore there are non-simply-connected leaves. The only
case of an $M$ that slithers around $S^1$ (or more generally, has
a codimension one foliation) such that all leaves are simply-connected
is when $M = T^3$. Note that in this case, the leaves of the
foliation are not hyperbolic.

For any non simply-connected leaf,
there must be a closed geodesic on the leaf---this follows from the
fact that the leaves have complete Riemannian metrics with
injectivity radius bounded from below. In a negatively curved metric,
any curve that is nearly geodesic is near a geodesic, so it is
easy to find a closed geodesic in each homotopy class, by curve-shortening.
Even in metrics for the leaves where there are no stipulations on
the curvature, we could apply a curve-shortening process to get
curves on leaves that are more and more nearly geodesic; such a curve
might not converge on its own leaf, but we could take a limit in $M^3$
to get a closed geodesic on some leaf.  This gives an invariant
measure.

For a non-singular flow such as $\gfl$,
we can factor any invariant measure as a transverse invariant measure
$\times\; dt$, where $dt$ denotes the measure of time along flow-lines
and coincides with arc length of a geodesic, in the case of the geodesic flow;
a transverse invariant measure is a measure locally defined on the
local space of orbits that is invariant under holomy.  The advantage
of converting to a transverse invariant measure is that it is 
a better topological invariant---intuitively, a transverse invariant
measure is something very similar to a closed orbit; it can be thought
of like a limit of $\R$-linear combinations of $1$-manifolds, and
is a special case of a closed
$1$-current, which is a version of a $\R$-1-cycle.  That is,
a transverse invariant measure gives a linear function on
$1$-forms on $T_1(\mathcal F)$, obtained by an iterated integral, first
integrating the $1$-form over flow-lines, then integrating with respect to
a transverse invariant measure.  Transverse invariance of the measure
is equivalent to the condition that this linear function vanishes on
$df$, for any function $f$.

Given a foliation with a Riemannian metric for its leaves,
let $\MG(\mathcal F)$ denote the space of transverse invariant measures
for the geodesic flow.  If the manifold is compact, this is the cone on
a compact convex set (using the weak topology on transverse invariant
measures).

Suppose now that we have a $3$-manifold $M^3$ that slithers around $S^1$,
with a Riemannian metric that is negatively curved on the leaves. 
Let $Z: \tilde M \to \tilde M$ be a homeomorphism that is a lift of
a generating deck transformation of $\R \to S^1$.  For any 
leaf $L$ of $\Tilde {\mathcal F}(s)$,
$Z(L)$ is a bounded distance from $L$.   Furthermore,
every leaf is sandwiched between two iterates $Z^k(L)$ and $Z^{k+1}(L)$.

\begin{lemma} \label{lemma: leaf hopping}
For each pair of leaves $L$ and $L'$ in $\tilde M$ and every
infinite geodesic $g$ on  $L$,
there is a unique geodesic $g'$ on $L'$ at a
a bounded distance from $g$.
\end{lemma}
\proof
Although the leaves are not quasi-isometrically embedded in $\tilde M$, 
they are properly embedded, since any fiber intersects any transverse
curve at most once.  In fact, the leaves are uniformly properly embedded,
in the sense that for any constant $A$ there is a constant $B$ such
that any two points on $L$ who have distance less than $A$ in $\tilde M$ have
distance less than $B$ along $L$. To see this fact, consider any sequence
of shortest geodesic arcs along leaves in $\tilde M$ whose
endpoints have distance in $\tilde M$ not
exceeding $A$. Adjusting by covering transformations and
passing to a subsequence, we may assume that the two endpoints converge. 
 Since the space of leaves in $\tilde M$ is Hausdorff, the pair of endpoints is
on a single leaf. The limit points are at a bounded distance on their leaves,
therefore, the length of the approximating arcs are bounded, and their lengths
converge to the distance between the limit points.

Let $C$ be the maximum separation between $L$ and $L'$ (so every point on
either leaf has a point within distance $C$ on the other).  Given any geodesic
$g$ on $L$, we can choose a sequence of points \set{p_i} on $g$ at uniform
intervals, say $1$, and find points $q_i$ within distance $C$ of $p_i$ on $L'$,
giving a path on $L'$ whose length is increased at most by some constant
factor. We can also go from geodesics on $L'$ to paths on $L$ that increase at
most by a constant factor.  This means that quasi-geodesics on either leaf
correspond to quasi-geodesics on the other.  In $\Hy^2$, or any other complete
metric of pinched negative curvature on $\R^2$, every quasi-geodesic is a
bounded distance from a unique geodesic.  This establishes the lemma.
\endproof
\begin{corollary}\label{corollary: circle identification}
If $s$ is a slithering of a 3-manifold $M^3$ around $S^1$ with
hyperbolic leaves, then there is a canonical
identification of the circles at infinity for all the leaves
of $\Tilde {\mathcal F}(s)$, giving a single $\pi_1(M)$-equivariant circle.
In other words, the circles at infinity for the leaves of $\mathcal F(s)$
forms a foliated circle bundle over $M$.
\end{corollary}
\proof
A geodesic on a leaf $L$ of $\tilde M$ is determined by a pair of
distinct points on its circle at infinity.   The bounded-distance
correspondence between geodesics on $L$ and geodesics on another leaf
$L'$ has to preserve this product structure, since two geodesics converge
to the same point at infinity if and only if their distance in
that direction stays bounded.
\endproof

Note that this circle bundle is not in
general isomorphic to the associated circle bundle of the slithering: for
example, in the case of a 3-manifold that fibers over $S^1$ with fibers of
negative Euler characteristic, the slithering bundle
is trivial, and the tangent circle bundle of the fibers is not.
However, section \ref{section: Anosov} discusses a situation when these
two bundles coincide.

\section{Canonical transverse flows} \label{section: pA}

A 3-manifold $M$ that has a slithering $s$
around $S^1$ has a foliation $\mathcal F(s)$,
but it also has a somewhat mysterious extra piece of dynamics,
the map $Z$ defined on the space of leaves of the foliation in $\tilde M$
which comes from the deck transformations of the universal cover of the
circle.  In this section, we will analyze the action of $Z$ on
geodesic currents on the leaves, enabling us to enhance $Z$
by constructing a connection for the slithering, canonical up to conjugacy by
a homeomorphism isotopic to the identity. In other words,
we will find a 1-dimensional foliation
transverse to $\mathcal F(s)$ and uniformly regulating it. Assuming
transverse orientability, if we isotope $M$ along  flow lines until
each point goes once around the circle, the result is a leaf-preserving
homeomorphism of $M$ that lifts to induce the automorphism $Z$ of
the space of leaves of the universal cover.

First, we can enhance $Z$ to define a map from $T_1(\mathcal F)$
to itself that gives an automorphism of the foliation by geodesics (flow-lines
of $\gfl$), using a standard trick.  In the unit tangent
space of the foliation of $\tilde M$, first construct
a $\pi_1(M)$-equivariant continuous map $f_1$
 that takes each point on a geodesic
on a leaf $L$ to some point at distance at most $A$ on the corresponding
geodesic on $Z(L)$.  This is readily constructed by making local
choices and using a partition of unity to average them. Now, define
$f(x)$ to be the average of $f_1(\gfl_t(x))$ where $t$ ranges over
a long interval $[-T, T]$. Since each flow-line has an affine structure,
averaging makes sense. It is easy to verify that $f_1$ is
quasi-monotone, that is, it eventually progresses in a single
direction. The quasi-monotonicity of $f_1$ implies monotonicity of
$f$.  (An alternative way to do this is to use triples of points on the circle
at infinity for the leaves.)

The automorphism $Z$ of $T_1(\mathcal F)$ induces an automorphism
of the space of transverse invariant measures, which we also denote
$Z: \MG(\mathcal F) \to \MG(\mathcal(F)$. If we denote $\PG(\mathcal F)$
the set of non-zero transverse invariant measures up to scaling,
then $Z$ acts as a projective endomorphism on $\PG(\mathcal F)$.

Every projective endomorphism of a compact convex set has at least
one fixed point, by Perron-Frobenius theory. In other words, it follows
from general principles that there is some transverse
invariant measure $\mu$ for $\gfl$ that is taken to a multiple of itself by $Z$.
We'll describe below a method to find such a $\mu$.
(A projective endomorphism of a compact convex set is
a generalization of a positive matrix, and a fixed point for the
transformation is a generalization of a positive eigenvector).  

Here's one fairly explicit way to obtain such a measure $\nu$. Let $\mu$ 
be any invariant measure for $\gfl$.
The \df{mass} $|\mu|$ of a transverse invariant measure $\mu$ is the
mass of $\mu$ converted back to an actual
measure $\mu \times ds$.  We consider all the images
$Z^k(\mu)$ for $k > 0$, and look at the maximal growth rate (or decay)
of its mass,
\[
g(\mu) = \limsup_{k \to \infty}
\left (\frac{ |Z^k(\mu)| }{ |\mu|}\right )^{1/k} .
\]
We will define a sequence of weighted averages
of $Z^l(\mu)$ that is more and more nearly a $1/g(\mu)$ eigenmeasure for $Z$.
We can do this by choosing weights for the first $2 N$ iterates such that
for the first $N$ terms
the weight of each term is slightly more than $1/g(\mu)$
times the weight of the preceding, for the last $N$ terms, slightly
less than $1/g(\mu)$ times the preceding. 
We can choose weights of this sort so that most of the total
comes from the middle $N$ terms of the sequence, provided $N$ is a number
so that for $l > N$ the estimates for $g(\mu)$ are not very much too
high, and for some of the middle terms they are near the limit. The
weighted sum, normalized to have mass $1$, is nearly a $1/g(\mu)$-eigenmeasure.
Any limit point of this sequence gives an eigenmeasure.

There is a kind of linking number between invariant measures for $\gfl$ that
will help give us a better geometric understanding of the action of $Z$.
This notion is a generalization of the intersection number of two measured
geodesic laminations or geodesic currents on a surface.
In the case of a compact surface $S$, the
\df{geometric intersection number} $i(\mu, \nu)$ is a symmetric bilinear
function of transverse invariant measures.
This geometric intersection number is the integral of the
product measure $\mu \times \nu$ over all intersection points of geodesics
on $S$.  If the injectivity radius of
$S$ is $\ge a$, then we can associate to any intersection point of two
geodesics $g, h$ the subset of $T_1(S) \times T_1(S)$ consisting of the
pair of segments of radius $a$ on $g$ and $h$; different intersections map
to disjoint sets.  Therefore, $a^2 i(\mu, \nu) < |\mu| |\nu|$. The
length of a geodesic is a special case of the intersection number with
the natural Riemannian volume element for the geodesic flow, that
is, the intersection number with a `random' geodesic. This fits into
a picture for Teichm"uller space analogous to the Lorentz
space model of $\Hy^n$, see \cite{Thurston:MinimalStretch}.

\begin{figure}[hbtp]
\includegraphics[width=.8\textwidth]{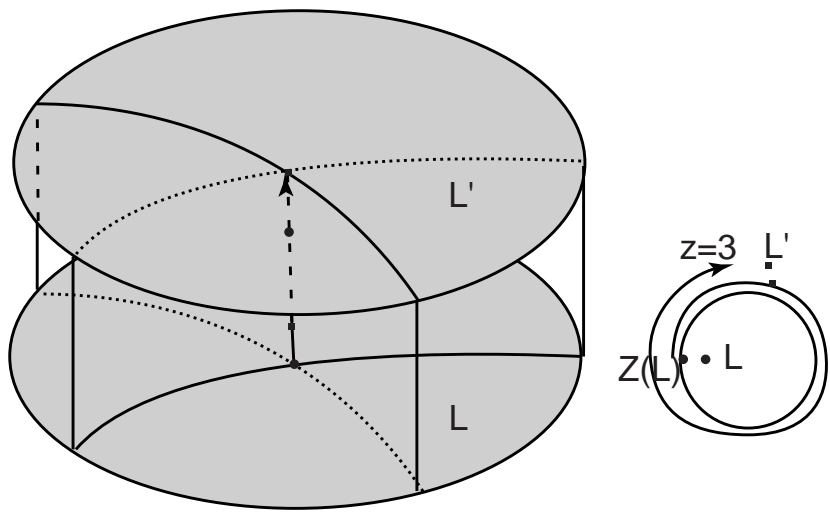}
\caption[Intersection of geodesics on
surfaces generalizes to crossing of geodesics in a slithering]
{ \label{figure: Lambda}
This is a sketch of two leaves in the universal cover of a slithering,
with two geodesics. The circles at infinity of the leaves are identified,
so we can tell whether or not the geodesics cross. Crossings modulo deck
transformations are identified by the homotopy class of a path
that joins them.  Given two geodesics currents $\mu$ and $\nu$
(measures on geodesics on the leaves),
we aggregate all the crossings into the linking
series $\Lambda(\mu,\nu)$ according to the $z$-value of a 
path that joins them.}
\end{figure}

To generalize intersection numbers
to the geodesic flow $\gfl$ for $\mathcal F$, we will
say that two geodesics on leaves in
$\tilde M$ \df{cross} if their projections to any single leaf intersect.
(See figure \ref{figure: Lambda}.)
We want to think about crossings modulo the action of $\pi_1(M)$.
The way to record this information (a pair of geodesics
modulo deck transformations) is with a homotopy class of
paths $\alpha$ between flow-lines of $\gfl$, where each endpoint of $\alpha$
can slide along but not leave its flow line.  A path of this sort
determines a crossing if for a lift to the universal cover,
the geodesic of its first endpoint crosses the geodesic
of its second.

For any crossing $\alpha$, we can give a measure of the height difference
of its endpoints, by setting $z(\alpha) = 2 k$ if $Z^k$ takes $\alpha(0)$,
lifted to the universal cover, to the leaf of $\alpha(1)$, and
$z(\alpha) = 2 k + 1$ if, lifted to $\tilde M$,
the leaf of $\alpha(1)$ is between the leaf
of $Z^k(\alpha(0))$ and $Z^{k+1}(\alpha(0))$. This definition makes
$z$ anti-symmetric in $\alpha$, that is, $z(\alpha) = -z(\alpha^{-1})$
where $\alpha^{-1}$ is the same path in the reverse direction.
Since the bounds for quasi-isometric comparisons between geodesics
on different leaves depend only on the value of $z$, every crossing
is represented by a path $\alpha$ whose length is
bounded as a function only of $z(\alpha)$.

Given two geodesics on leaves of $\tilde M$, we can
quantify crossing data by grouping 
pairs of geodesics according to the value of $z$.
We encode this information in a formal Laurent series, integrating
over all crossings $\alpha$:
\begin{equation}
\Lambda( \mu, \nu) = \int_{\alpha}  t^{z(\alpha)}
\mu \times \nu
\end{equation}
Convergence of the integral can be checked similarly to convergence
in the definition of intersection number of laminations on surfaces.
If the injectivity radius of $M$ is $a$, then given an arc $\alpha$,
any other arc joining points on segments of radius $a$
on the geodesics of its endpoints and staying within $a$ of $\alpha$
represents the same crossing. The space of arcs of length less than
$B$ is compact; each coefficient is dominated by an integral of a 
locally bounded measure over a compact space, therefore is bounded.
This reasoning shows that for each $k$ there is some constant $C_k$
such that the coefficient of $t^k$ in $\Lambda(\mu, \nu)$ is
not greater than $C_k \left | \mu \right | \left | \nu \right |$.

The constant term of $\Lambda(\mu, \nu)$ measures the actual
intersections.
\[
\Lambda(\mu, \nu)(t) = \Lambda(\mu, \nu(1/t).
\]
The coefficient of $t$ measures how much $\mu$ intersects
$\nu$ when each geodesic is swept through $M$ to its image under $Z$.
Other coefficients are similar.  This linking series is not continuous
as a function of $\mu$ and $\nu$: when even terms are non-zero, it's
sometimes possible for $\mu$ and $\nu$ to have arbitrarily
small perturbations where some or all of the weights on even
terms jump to neighboring odd terms.  

If $\mu$ and $\nu$ are invariant measures, note that
\[ \label{equ:Z Lambda}
\Lambda(\mu, \nu) = t^2 \Lambda(Z(\mu), \nu) =  t^{-2} \Lambda(\mu, Z(\nu)) .
\]
This implies that the coefficients of $\Lambda(\mu, \nu)$ grows at
most exponentially fast.  In other words, the portion with
positive exponent converges in some neighborhood of $t=0$, and the
portion with negative exponent converges in a neighborhood of $t=\infty$.
Of more immediate significance is the consequence that
the linking series is invariant when both arguments are transformed by
$Z$.  This implies that if $\mu$ is any eigenmeasure
for $Z$ with eigenvalue not 1, then
$\Lambda(\mu, \mu)$ is identically $0$.  It also implies that if
$\mu$ and $\nu$ are any two eigenmeasures for $Z$ such that
$\Lambda(\mu, \nu) \ne 0$, then their eigenvalues are reciprocal.

\begin{proposition} \label{proposition: bounded or simple}
Either
\begin{itemize}
\item[A.]
The action of $Z$ on geodesic currents is globally bounded, that is,
there is some constant $K$ such that for every transverse invariant measure
$\mu$ for $\gfl$ and for every integer $k$, 
\[
\left |Z^k(\mu)\right | < K \left |\mu\right |,
\]
or
\item[B.]
there is a transverse invariant measure $\nu$ for $\gfl$ that is an
eigenmeasure of $Z$ and such that
$\Lambda(\nu, \nu) = 0$
\end{itemize}
\end{proposition}
\proof
Suppose that the action of $Z$ on geodesic currents is not globally
bounded, so there is a sequence \set{\mu_i} of transverse invariant
measures for $\gfl$ of mass $1$, and integers $\set{k_i}$ such that
\set{\left | Z^{k_i} (\mu_i) \right |} tends to $\infty$.
Then the sequence of normalized transforms
\set{\mu'_i = Z^{k_i}(\mu_i) / \left | Z^{k_i} (\mu_i) \right |}
has the property that all the coefficients
of $\Lambda(\mu'_i, \mu'_i)$ tend to $0$ as $i$ tends to $\infty$,
in light of equation \ref{equ:Z Lambda}.  We can pass to the limit,
and obtain an invariant measure $\mu$ for $\gfl$. Even though $\Lambda$
is not in general continuous, it follows in this case that
$\Lambda(\mu, \mu) = 0$, since the terms of the limit only have
contributions from the limits of neighboring terms.  In other words,
there are no crossings at all between
geodesics in the support of $\mu$, in no
matter what homotopy class $\alpha$.  The set $\sigma$
consisting of all geodesics
on leaves of $\mathcal F$ that do not cross geodesics in the
support of $\mu$ forms a closed set, invariant under $\gfl$ and under
$Z$. Therefore, there is a transverse invariant measure $\nu$ for $\sigma$
that is projectively invariant by $Z$.
\endproof

\begin{remark} \label{remark: 1 but grows}
It can happen that the masses of images of invariant measures are unbounded,
yet the only positive eigenmeasures for $Z$ have eigenvalue $1$. This
happens, for example, with the mapping torus of a Dehn twist $T_\gamma$
around a non-trivial curve $\gamma$ on a negatively curved surface.
In this case, $Z$ acts on fibers as $T_\gamma$; it represents isotoping
of the mapping torus once around the circle.  For any curve $\beta$
that crosses $\gamma$, the
images $Z^k(\beta)$ grow to infinity in length, while
the normalized transverse invariant measure supported on
$Z^k\beta$ converges to $\gamma$, which is an eigenmeasure of
eigenvalue $1$.  
\end{remark}

A slithering $s$ of $M^3$ is \df{reducible} when there is an embedded
incompressible torus or Klein bottle where $s$ induces a slithering.
A torus of this form is a \df{reducing torus} (or Klein bottle if it's
a Klein bottle.)  It is known from foliation theory
that in a foliation of a 3-manifold without Reeb components,
any incompressible  torus is isotopic to a transverse torus,
or a leaf in the special case of a foliation
whose leaves are all tori. Consideration of the various cases shows
that a torus transverse to $\mathcal{F}(s)$ is a reducing
torus, unless $M$ fibers over $S^1$ or $I$
with fiber a torus (which does not
imply that $s$ is itself a fibration.)

Before proceeding further with the generic case B of proposition
\ref{proposition: bounded or simple}, the non-generic case A
has an interesting story:

\begin{theorem} \label{theorem: convergence groups}
\textbf{(Convergence group theorem, Gabai \cite{Gabai:Convergence} and
Casson and Jungreis, \cite{Casson:Jungreis:Convergence})}
In case A of proposition \ref{proposition: bounded or simple},
$M$ is a Seifert fiber space and the slithering is a foliated circle
bundle.
\end{theorem}
\proof
\emph{Note:}  this is not a new proof, but simply a reduction of the current
statement to a standard established form. 

The analysis of this case is a consequence of a long, winding
series of developments starting with Nielsen, who settled the case
mapping torus case, which is equivalent to---Nielsen
showed that if some power a diffeomorphism $\phi$ of a surface is homotopic
to the identity, then $\phi$ is isotopic to a diffeomorphism of finite
order, which implies that the mapping torus of $\phi$ can also be
described as a Seifert fiber space.  We will not attempt to recount
the history, nor to fit the present application naturally into context.
However, it's worth remarking that the hypothesis (A) is close to
the condition that $\pi_1(M)$ acts as a group of uniform quasi-isometries
of the hyperbolic plane, which is close to the condition that it
acts on the circle as a convergence group, meaning its action
on triples of points on the circle is properly discontinuous. It's also
worth remarking that the set of triples of points on a circle,
modulo a convergence group, comes naturally equipped with a slithering
around $S^1$.

We will content ourselves here with deriving the theorem logically
from established knowledge. For this, it will be sufficient to find
an immersed incompressible torus.  Assuming $M$ is not $T^3$
(where the theorem is more trivially true), not all the leaves of
$\mathcal F(s)$ are simply connected. Let $\gamma$ be a non-trivial
curve on any leaf. Choose a base-point $* \in M$, and for each image
$Z(\gamma)$ connect $\gamma$ by a shortest path to $*$.  This
gives a sequence of homotopy classes of bounded length, so they
repeat infinitely often.
We can form a long immersed cylinder connecting all the images
$Z^k(\gamma)$ and
intersecting each leaf in a geodesic, not necessarily closed, but
converging to the same end points on $S_\infty^1$.
The homotopy information tells us that this
this cylinder eventually joins up, forming a torus or Klein bottle.

If $\mathcal F(s)$ is reducible, then it decomposes into geometric
pieces by the geometrization theorem for Haken manifolds. Hyperbolic
pieces are incompatible with immersed incompressible non-boundary-parallel
tori.  Taut foliations of Seifert
fiber spaces have been understood for some time---the Haken cases were
analyzed in \cite{Thurston:thesis}, and the non-Haken cases followed from
later developments. In general, a taut foliation of a Seifert fiber is
isotopic to be transverse to the fibers provided
the fiber is not homotopic to a leaf. If the base is a negatively-curved
orbifold, the only possibility for $\mathcal F(s)$ to be isotopic so
that it is transverse to the fibers. It can happen that the leaves
of $s$ are `vertical' when the base is Euclidean,
but then (given condition A) there is some other Seifert fibration
transverse to the leaves.

If $s$ arises by sewing together Seifert fiber pieces along tori
in a way that fibers cannot be chosen to align, then a geodesic
current that crosses an offending torus violates hypothesis A (just
as in remark \ref{remark: 1 but grows}.)

If $M$ is torus-irreducible but has an immersed incompressible torus,
the culmination of the long development mentioned above implies
that $M$ must be a Seifert-fiber space.
\endproof

We need to think about the topology as well as the measure
theory of the action of $\pi_1(M)$ on $\gfl$.  To develop the
topological picture a little further, consider any
$L$ in $\tilde M$. The fundamental group of $M$ acts
on the set of geodesics in $L$.   Consider any closed invariant
subset $I$ for this action.  Form an $\epsilon$-neighborhood $I_\epsilon$
(on $L$) of the union of geodesics in $I$. We claim that
if $I_\epsilon$, or any
connected component of $I_\epsilon$, does not limit to all of $S_\infty^1$,
then $s$ is reducible.  To see this, let $X \subset S_\infty^1$ be any
proper subset that is the limit set for some component of $I_\epsilon$,
and form the convex hull $H$ of $X$ in the hyperbolic metric on $L$. Each
element of $\pi_1(M)$ takes $H$ to itself or to a set that is
disjoint except possibly on one boundary line. More specifically,
it is possible to have two distinct boundary lines in a single
$\epsilon$ disk on $L$, but not three.  Now consider all the
images of boundary lines of $H$ back in $M$, transported to all leaves of
$\mathcal F(s)$.  The union is compact, since $M$ is compact and each sheet of
the surface swept out by the  boundary lines has a uniform neighborhood
intersecting at most one other sheet.  Therefore, the union is a compact
surface. The surface has an induced slithering, which means its
Euler characteristic is $0$, so it is a torus or Klein bottle.

Note in particular
that if there is some $\epsilon$ such that $I_\epsilon$ is not
connected, then any component of $I_\epsilon$ has a limit set  $X$
that is a proper subset of $S_\infty^1$, so $s$ is reducible.

We will say that a $\pi_1(M)$-invariant set $I$ of geodesics 
\df{fills} $M^3$ if for every $\epsilon$, $I_\epsilon$ is connected,
and its limit set is the entire circle at infinity. As we have seen,
if $s$ is irreducible, then every $I$ fills $M$.

If there are crossings in $M$ among the geodesics in $I$, then there
may not be a clear choice of a canonical form for the
immersed surfaces swept out by $I$ in $M$---this is the usual problem,
that whenever three
or more lines cross each other, there are multiple patterns in which
they can cross, and it is hard to choose among them. However, when
there are no crossings, $I$ sweeps out a topologically well-defined
2-dimensional lamination $S(I)$ in $M$, intersecting each leaf in
the geodesic lamination covered by $I$. 
\begin{lemma} \label{lemma: gaps exist}
If $I$ is a closed $\pi_1(M)$-invariant subset of geodesics
with no crossings, then the gaps of $I$ are dense, that
is, no neighborhood in $\Hy^2$ is foliated by geodesics in $I$.
\end{lemma}
\proof
If any open set is swept out by geodesics, one or the other
endpoint of the geodesics actually moves.
If we go toward a non-constant endpoint in $\Hy^2$, the geometric limit
is a foliation of $\Hy^2$ by geodesics having one endpoint constant
and the other endpoint completely traversing the remainder
of the circle. Since $M$ is compact, this behavior would actually
occur in the intersection of $S(I)$ with some leaf, and hence
in every leaf. It follows that $\pi_1(M)$
would have a fixed point on $S_\infty^1$, so the action of $\pi_1(M)$
on a leaf is solvable. There is enough information in this
picture to determine that $M$ would have to be commensurable with
the mapping torus of an Anosov diffeomorphism of the $T^2$, as developed
in the analysis of Dehn surgeries on the figure eight knot in
Thurston, \cite{Thurston:GT3M}.  However, 
the foliations compatible with this foliation
are not of the form $\mathcal F(s)$.
\endproof

It is easy to see that $S(I)$ is an essential lamination.
Each gap in the lamination
$I$ sweeps out a submanifold of $M$, which is a 
gap in $S(I)$ .  If we truncate the 3-dimensional gaps where the
geodesic sides approach within $\epsilon$, we obtain a finite
number of compact 3-manifold with walls that alternate being
geodesics on their leaves and $\epsilon$-hops from one geodesic
wall to another.  If $s$ is irreducible, then the gaps
are solid tori, in the form of ideal polygons that wrap around
through $M$ and return with some rotation.

\begin{proposition} \label{proposition: lamination regulates}
If $s$ is an irreducible slithering that has an invariant set $I$
of non-crossing geodesics, then $s$ is regulated by a flow $\phi_t$
that preserves the 2-dimensional lamination $S(I)$.
\end{proposition}
\begin{proof}
We can first construct the flow $\phi$ on $S(I)$ by first making a rough
but bounded guess, then averaging along the geodesics.
If the flow is chosen continuously on $S(I)$, we can easily extend continuously
to the solid torus gaps.
\end{proof}

\begin{proposition} \label{proposition: reciprocal eigenvalues}
Let $s$ be an irreducible slithering of a closed three-manifold $M^3
$ around $S^1$ which is not in case A of
\ref{proposition: bounded or simple} (\emph{i.e.}, is not a Seifert
fiber space, by \ref{theorem: convergence groups}).
 Let $\lambda$ be the largest sustained growth factor
for any lamination under positive or negative iterates of $Z$, that is,
\[
\lambda = \limsup_{|k| \to \infty}
\left ( \frac{ |Z^k(\mu)|}{|\mu|} \right )^{1/k}
\]
Then there are two transverse invariant measures
 $\mu_s$ and $\mu_u$ for $\gfl$ that are
$\lambda$ and $1/\lambda$ eigenmeasures for $Z$.
\end{proposition}
\proof
We can use the procedure described above to find an eigenmeasure whose
transverse invariant measure shrinks in a direction that the growth factor
is attained.  Say this is $\mu_u$, whose transverse invariant measure
grows for negative $k$.

Now let $\nu$ be any measure such that $\Lambda(\mu_s,  \nu) \ne 0$,
from a construction of
section \ref{section: geodesic currents}.
Because $\mu_s$ is a $\lambda^{-1}$ eigenmeasure of $Z$, 
this series has the form $(a_0 + a_1 t) \sum_k \lambda^{-k} t^{2k}$.
In other words, $\nu$ crosses $\mu_s$ more and more for negative time,
with the measure of crossings growing by a factor of $\lambda$ at each
stage.  This can only happen if the mass of $\nu$ grows geometrically,
comparable to $\lambda^{-k}$ for $k < 0$.  We can construct a
$\lambda$-eigenmeasure by taking any limit point $\mu_s$
of weighted combinations of these
images under $Z^{-k}$.   
\endproof

So far, we have not logically deduced 
the geometric relationship between $\mu_u$ and $\mu_s$. It
is logically consistent with proposition
\ref{proposition: reciprocal eigenvalues} and its proof
(even if this may seem bizarre geometrically)
that $\mu_u$ and $\mu_s$ are mutually singular measures
that are physically supported on the same invariant
subset of $\gfl$.  Our next task is to analyze this geometry,
so that, in particular, we will be able to use eigenmeasures for $Z$
to construct a pseudo-Anosov flow transverse to $\mathcal F(s)$.

We suppose that $s$ is an irreducible slithering of a closed 3-manifold
$M$ with transversely oriented foliation $\mathcal F(s)$,
and that $\nu$ is a transverse invariant measure for $\gfl$ and
eigenmeasure for $Z$ with eigenvalue $\lambda < 1$. Let $I_\nu \subset L$
be the support of $\nu$.  Define $M_{L+}$ to be closed halfspace of
$\tilde M$, on the positive side of $L$,  and $M_{L-}$ similarly.
If we project geodesics in $M_{L+}$ to $L$, this induces a
map on transverse invariant measures; since $\lambda < 1$, the contributions
of $\nu$ to the projection decrease exponentially as a function of the
height function $z$, so this gives a well-defined transverse invariant
measure for $I_\nu$. Let's call this image $L_+(\nu)$.
(We could also have projected just
the contribution from a fundamental domain $0 \le z < 2$, which would give
a transverse invariant measure for $I_\nu$ no matter
what the value of $\lambda$.)

Note that the corresponding measure on $Z(L)$ is the measure on $L$ 
multiplied by $\lambda$:
\[
Z^k(L)_+(\nu) = \lambda^k L_+(\nu).
\]

\begin{lemma} \label{lemma: no atoms}
$L_+(\nu)$ has no atoms, that is, there is no single geodesic of $I_\nu$
with positive transverse measure.
\end{lemma}
\proof
If there were any atoms for $L_+(\nu)$,
its images under negative powers of $Z$ would have arbitrarily
large mass; this is incompatible with compactness of $M$ and the
finite mass of $\nu$.
\endproof

Let's focus on one gap $G$ on $L$, and look at what
happens along one of its sides $g$ (where $g$ is a geodesic).  Choose
a short arc $J$ transverse to $I_\nu$
and connecting $G$ to another
gap. Then $J \cap I_\nu$ is a Cantor set, since gaps are dense and
there are no atoms to $L_+(\nu)$.
Define $\mathcal G \subset \R$ to
be the collection of values of the $\nu$-measures along $J$ from $G$ to
other gaps; then $\mathcal G$ is a dense subset of the interval
$[0,\nu_J]$.

There are infinitely many gaps in $I_\nu$, but each of them is
an ideal polygon that is a lift of the intersection of one 
of the finitely many solid tori gaps in $S(I_\nu)$.  This implies that
there is bounded variability of the geometry and of the transverse
invariant measure among all the gaps of $I_\nu$ wherever they
appear, among all the leaves of $\tilde M$.

Let $S(G)$ be the solid torus gap in $M$ swept out by $G$. We can
follow the solid torus once around, giving a return map of $G$ to itself.
The return map $R_G$ of $G$ might take $g$ to a different side of $G$, but
some iterate $R_g = R_G^p$ takes $g$ back to itself. The return map
is not in general a power of $Z$, and will not in general take
$L_+(\nu)$ projectively to itself, but nonetheless the return map is sandwiched
between two powers of $Z$: if
$\alpha$ is an arc in the homotopy class of $R_g$, which can
be identified with a deck transformation of $\tilde M \to M$ that does
the right thing to $g$, then it is sandwiched between
$\floor{z(\alpha)/2}$ and $\ceiling{z(\alpha)/2}$ power of $Z$, translating to
the inequalities
\[
\lambda^{\ceiling{z(\alpha)}} L_+(\nu) \le
R_g^* L_+(\nu) \le
\lambda^{\floor{z(\alpha)}} L_+(\nu).
\]
Since the transverse measure for $I_\nu$ shrinks when
pushed forward by the return map, this says that the leaves of $I_\nu$
have to spread further from each other. We are aiming to construct
particular return maps that balance the separation of
geodesics by quasi-uniformly shortening them.

We may assume that $J$ is short enough so that it cuts each gap 
it intersects other than $G$ into a single tip
of the ideal polygon on one side, and
the thick part of the ideal polygon on the other.
$R_g$ sends $J$ to some arc $R_g(J)$ transverse to $I_\nu$; no matter
how $J$ was chosen,
at least some initial segment of this image arc crosses the same leaves
as an initial segment of $J$.
We may replace $J$ by a segment that has an
$I_\nu$-preserving isotopy to an initial segment of $R_g(J)$.

\begin{figure}[thbp]
\includegraphics[width=.9\textwidth]{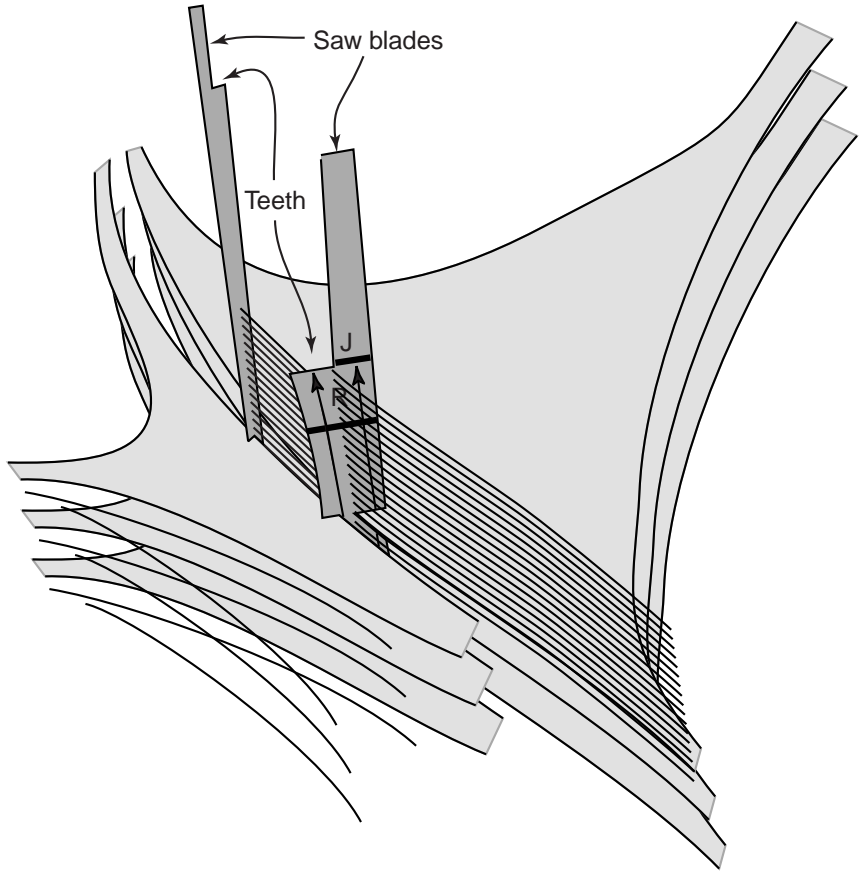}
\caption{ \label{figure: saws}
Saws cutting through laminations}
\end{figure}

Now construct an annulus $S(J)$
transverse to $S(I_\nu)$ by sweeping $J$ around
$S(G)$, always intersecting the same set of leaves until it
returns via $R_g$; at that point, glue an initial segment to $J$.
The figure formed is like a one-tooth saw-blade slicing through 
layers of $S(I_\nu)$ (figure \ref{figure: saws}.)

We can continue in the same vein, to
construct a family $S(J_i)$ of disjoint saw-blade annuli, one
for each annular face of each solid torus
gap of $S(I_\nu)$. Only the tooth of a saw-blade cuts through $S(I_\nu)$;
except for the tooth, the rim is in a gap.  

Every leaf of $S(I_\nu)$ is dense in $S(I_\nu)$---otherwise, its
closure would either contain all of $S(I_\nu)$ except for isolated leaves,
which we have ruled out by lemma \ref{lemma: no atoms}, or it would
be small enough to give us a reducing torus.
It follows that the saw blades $S(J_i)$ have sliced each geodesic of $I_\nu$
on each leaf of $\mathcal F(s)$ into bounded intervals. 

We can group these intervals of geodesics
according to their homotopy class \emph{rel}
saw blades. On any one leaf of $\mathcal F(s)$,
they group into a locally finite collection
of parallel bundles. Topologically, if the parallel arcs in one bundle
are squeezed together to one arc, the collection of arcs plus
intersections of saw-blades divides the surface into compact simply-connected
regions, which become polygons if the saw-blade intersections are collapsed
to points.  The polygons come from gaps of $I_\nu$ and from ends of
the saw-blade intersections (or both).

In three dimensions, a set of parallel arcs sweeps out sheets (like a
layer pastry). In the downward direction, sheets can
run into saw teeth, where they are cut into two pieces.  If there were
any family of arcs that could be isotoped downward forever, then the
arcs would have to stay bounded in length forever.  The arc defines
a homotopy class of paths between the cores of the two solid tori;
there are only finitely many homotopy classes of bounded length, so
eventually this homotopy class repeats, joining to form an annulus
that gives a homotopy of some power of the core of one
solid torus to some power of the core of some solid torus (possibly the
same). The intersection of the annulus with leaves of $\mathcal F(s)$ 
gives homotopy classes of arcs between two gaps; however,
because $\lambda \ne 1$, these arcs would have zero intersection number
with $L_+(\nu)$, which is impossible.  Therefore, every sheet of arcs
is split by a saw-tooth in the downward direction.
In the upward direction, every sheet that is not an original gap boundary
eventually runs off the edge of a saw tooth, where it merges with other sheets.
In the three-dimensional picture in $M$, only finitely
many homotopy classes of arcs occur, where each homotopy class $\beta$
serves as an index pointing to a family $P_\beta$
of parallel rectangular sheets.  Let $t_\beta \in [0,1]$ be a parameter
for the vertical direction of $P_\beta$, so that $t_\beta$ locally
parameterizes the leaves of $\mathcal F(s)$ within a rectangular
solid that encloses $P_\beta$.

\medskip
At this stage,  we can apply the automorphism $Z$, that is,
turn on the saws, after first adjusting the blades into proper alignment.
For each solid torus gap, attach the hubs of all its saw blades to
the core circle.  The image of a saw-blade under $Z^{-1}$ intersects
fewer leaves; it is isotopic to a subset of itself. If we modify
$Z$ using such an isotopy,
then the image of the each saw blade under $Z$ contains the saw blade.
We can also arrange the flow $\phi$ to be tangent to the
saw blades and to their edges, so that $Z$ takes each saw blade to the
same flow line. Under iteration of $Z$, the blades expand indefinitely,
by consideration of the measure $L_+(\nu)$.  Their edges remain
bounded in length, since they traverse a bounded $z$-length.  As
the saw blades expand, they remain transverse to the intervals of $I_\nu$,
chopping them successively into more and more pieces.  In other words,
the blades appear to converge to a lamination that has a finite number
of leaves with a singularity at the core of each solid torus gap of
$S(I_\nu)$.  To make a precise definition
of a limiting lamination, we can consider the order type of intersections 
of the saw blades with the arcs $\beta$.
The union of all intersections gives a
countable linearly-ordered set on each homotopy class of
arcs. The order completion
of any countable linearly-ordered set
(completion with Dedekind cuts, adjoining an upper and lower endpoint
if needed)
is a compact linearly-ordered set $X_\beta$
that has an order-preserving embedding in the unit interval unique up
to homeomorphism (and might well be an interval.)  For each
$\beta$, form the product of $X_\beta$ with a rectangle; we'll assemble
these pieces to form a singular foliation transverse to $\mathcal F(s)$
and to $S(I_\nu)$.  For each $\beta$, we have a rectangular solid, with
a foliation or lamination already in two coordinate directions. The
new leaves have automatic gluing maps at their top and bottom faces,
since these come from the cutting edges of the saw blades, where each
intersection continues through.  The horizontal faces work similarly.

For the sake of symmetry, we can blow up the resulting lamination along
its singular leaf, replacing the singular leaf with a solid torus gap,
having new leaves in the form of annuli for each `side' of the singular 
leaf. We'll call this new lamination $S^*(I_\nu)$.
The intersections of leaves of $S^*(I_\nu)$ with leaves of $\mathcal F(s)$
form a 1-dimensional foliation $I_\nu^*$, that is, $S^*(I_\nu) = S(I_\nu^*)$.
Each leaf of $I_\nu^*$, lifted to $\tilde M$,
converges to distinct points at its two endpoints, on account of its
transversality to the geodesic lamination $I_\nu$. (The leaves of $I_\nu$
describe neighborhood bases for all points at infinity except the vertices
of its gaps; the leaves of $I_\nu^*$ do not stay in any gap of $I_\nu$,
and cross $I_\nu$ at angles bounded from below.)  Therefore, each leaf
of $I_\nu^*$ is homotopic to a unique geodesic. The homotopies can be
made uniformly bounded, using the compactness of $M$. We straighten
out all the leaves of $I_\nu^*$ to geodesics; it doesn't matter if
conceivably this collapses multiple leaves to a single leaf, since
the process is invariant from leaf to leaf.  The resulting lamination
is still transverse to  $I_\nu$, so we can simply make this adjustment,
keeping the same names $I_\nu^*$ and $S^*(I_\nu)$ for the 1-dimensional
and 2-dimensional laminations.

Now we apply previous constructions, to obtain
an invariant measure $\nu^*$ on $I^*_\nu$ which is a $1/\lambda$
eigenmeasure for $Z$. We obtain a measure $L_-(I_\nu^*)$
on each leaf.
By examining the possible gaps, which necessarily are ideal polygons,
it's clear that its
support must be all of $I^*_\nu$.

To summarize and slightly extend:

\begin{theorem} \label{theorem: transverse pseudo-Anosov}
Let $s: \Tilde M \to S^1$ be a slithering of a compact $3$-manifold $M$.
Then $M$ splits along some (possibly empty)
family of reducing tori into a finite number of pieces.  On each
of the pieces, the induced slithering is either a foliation transverse
to the fibers of a Seifert fiber space, or it admits a transverse
pseudo-Anosov flow $\phi_t$ (or possibly a
pseudo-Anosov line field, if $\mathcal F(s)$ is not transversely orientable)
whose stable and unstable foliations are uniquely determined by $s$.
\end{theorem}
\proof
In the case of a closed manifold with an irreducible slithering, we
are nearly done. We have already addressed the case that lengths
of invariant measures for the geodesic flow transform boundedly.

We should recall here that
if the leaves of $\mathcal F(s)$ do not have hyperbolic type, then 
$\mathcal F(s)$ admits a transverse invariant measure (for instance,
by Candel, \cite{Candel:Uniformization}), and $s$ is
a perturbation of a mapping torus of a diffeomorphism of $T^2$
(or possibly a bundle over $S^1$ with fiber a Klein bottle, or a
bundle over the $I$ orbifold, if we do not make assumptions
of orientability and transverse orientability).
In this situation, the theorem is easily verified.

If lengths of invariant measures transform unboundedly but have eigenvalue $1$,
we still get an eigenmeasure $\nu$ for $Z$ whose support has no crossings
(proposition \ref{proposition: bounded or simple}.)
The eigenvalue $1$ implies there is a transverse invariant measure for
$S(I_\nu)$.  Such a measure is expressible as a cohomology class with
coefficients twisted by the transverse orientation of $S(I_\nu)$. It can be
perturbed to a rational class, which gives a perturbed 2-dimensional lamination
having nearby support, all of whose leaves are closed.  These leaves have
foliations induced from $\mathcal F(s)$, so they are toruses, and $s$ is
reducible.

If there is an eigenmeasure $\nu$ for $Z$ of eigenvalue $\ne 1$, we have
the two laminations $S(I_\nu)$ and $S^*(I_\nu)$ constructed in
the preceding discussion.
By a familiar technique (see e.g. \cite{Thurston:GT3M}),
given $\mathcal F(s)$, $S(I_\nu)$ and $S^*(I_\nu)$, we obtain
a decomposition of $M$ consisting of intersections of leaves and/or
gaps from the three. If we collapse the elements of this decomposition,
we obtain a homeomorphic manifold where the laminations are singular
foliations, intersecting in $1$-dimensional foliation. Every arc
of the flow of $z=1$ multiplies $L_+(\nu)$ by $\lambda$ and $L_-(\nu*)$
by $1/\lambda$. 
We can adjust the transverse metric so that
the transverse flow steadily compress in one normal direction while 
expanding in the other, making it a pseudo-Anosov flow.
Given uniqueness of the pseudo-Anosov foliations,
the non-transversely-orientable case
can be obtained from the transversely oriented double cover.

Any other invariant measure $\mu$ for the geodesic flow has to have
crossings either with $\nu$ or with $\nu^*$, so its eigenvalue is
determined, and it can have no crossings with the other of the two.
This implies uniqueness of the geodesic laminations
$I_\nu$ and $I_\nu^*$, which implies uniqueness of the pseudo-Anosov
foliations.

\medskip
When $s$ has torus boundary, for each boundary component, we can take
a horoball in $\Hy^3$ modulo $\Z^2$, glue it onto the given boundary component,
and extend the foliation by coning to the cusp. The point is
that when we uniformize the resulting leaves, we will obtain well-behaved
metrics that behave like coverings of complete
hyperbolic surfaces of finite area.  We obtain an identification
of the circles at infinity of leaves in the universal cover---we 
actually have the best control of the quasi-constants inside cusps.
This allows us to make all the constructions, projecting geodesics
from one leaf to another, and in particular defining the action of $Z$ on the
geodesic flow $\gfl$ of the leaves.

As before, if there is an eigenmeasure for $Z$ with eigenvalue $1$,
the slithering is either reducible, or $M$ is a Seifert fiber space.

When there is an eigenmeasure $\nu$ of eigenvalue $\ne 1$, then $\nu$
can have no crossings with itself.  
This forces each end of each geodesic in the support of $\nu$ either to
avoid a neighborhood of the cusps, or to head straight for the cusp:
any geodesic that goes very far toward
the cusp without going all the way wraps around and crosses itself.
By Poincar\'e recurrence, almost all measure leaving a cusp ends
up going back into a cusp. Since there are countably many cusps, this
means almost all measure near cusps is supported on atoms. But
there can be no atoms since $\lambda \ne 1$, therefore $\nu$ has
compact, bounded support.

When we now look at the two-dimensional lamination $S(I_\nu)$, there is
an additional kind of gap, enclosing a cusp.  These have a strong
resemblance to the other, solid torus, gaps:
they have the form of the suspension
of a map of a punctured ideal
polygon to itself, where the punctured ideal polygon has one or more sides.
We can use these to build a flow $\phi_t$ that regulates
$s$ and is tangent to $S(I_\nu)$.

We can construct saw blades for the cusp gaps, just as for the solid torus
gaps.  We cone the hub of each cusp saw blades to the cusp (attaching
a pseudosphere) before turning the saws on.
The construction of $S^*(I_\nu)$ goes through just as before.
\endproof

\section{Peano curves}
\label{section: Peano curves}

Suppose $\mathcal F$ is a taut foliation
of a hyperbolic 3-manifold $M^3$.  It is known (Fenley,
\cite{Fenley:quasi-isometric}) that not all of the leaves 
can be quasi-isometrically embedded:
there exist geodesic paths in any Riemannian metrics for the leaves
that can be shortened in $M^3$ by arbitrarily large factors using homotopy
in $M^3$ \emph{rel} endpoints.  In other words, the geometry of the
leaves is very far from the geometry of a hyperbolic plane in hyperbolic
3-space---so, what do they look like? 
In particular, in $\Tilde M$, do the leaves extend continuously
to give a map of $\Hy^2 \cup S_\infty^1$ to $\Hy^3 \cup S_\infty^2$,
and if so, what is the topology and geometry of this map?

The geometry is already interesting in the simplest taut foliations
of hyperbolic 3-manifolds, the foliations by the fibers for three-manifolds
that fiber over $S^1$.  This case was resolved
by Cannon and Thurston (\cite{Cannon:Thurston:peano}),
where it was shown that the universal coverings of fibers extend to
define sphere-filling `Peano' curves. Fenley (\cite{MR93k:57030})
generalized this result to the case of depth one foliations, that is,
foliations such that every leaf is either closed, or accumulates
only on closed leaves. Fenley showed in the depth one case that all
leaves converge at infinity, but the limits are
not sphere-filling curves except when the closed  leaves
are fibers of a fibration over $S^1$. The `typical' behavior at infinity
is for the depth zero (closed) leaves to limit as circles, and the
depth one leaves to limit to curves whose image is a swiss cheese,
but there are also other possibilities, depending on the nature
of the characteristic $I$-bundles for $M^3$ split 
along various leaves.

We will show here that the behavior of slitherings at
$S_\infty^2$ is like the case of manifolds that fiber over $S^1$.

\begin{theorem} \label{theorem: Peano curves}
Let $s: \Tilde M \to S^1$ be a slithering of a compact 3-manifold $M$,
whose interior has a complete hyperbolic metric of finite volume.
\begin{enumerate}
\item[a.]
The universal covers of the leaves of $\mathcal F(s)$, lifted
to $\Hy^3$, extend to give continuous maps
\[ \Hy^2 \cup S_\infty^1 \to \Hy^3 \cup S_\infty^2 . \]
These maps respect the identification (corollary
\ref{corollary: circle identification}) of their circles at infinity.
\item[b.]
The universal covers of
leaves of the stable and unstable laminations associated with the
transverse pseudo-Anosov flow
(theorem \ref{theorem: transverse pseudo-Anosov}), lifted to $\Hy^3$, extend
to give another set of continuous maps
\[ \Hy^2 \cup S_\infty^1 \to \Hy^3 \cup S_\infty^2 . \]
\item[c.]
If $M$ is closed, then the universal coverings of leaves of
the stable and unstable pseudo-Anosov laminations
are quasi-isometrically embedded in $\Hy^3$. 
\item[d.]
If $\boundary M$ is non-empty (necessarily it consists of
tori and Klein bottles) then the universal coverings of
leaves of the stable and unstable
laminations are not quasi-isometrically embedded in $\Hy^3$. For
any stable or unstable leaf $l \subset M$ having a non-trivial
closed loop homotopic to the boundary, the 
two endpoints of the universal cover of the loop are identified
at $S_\infty^2$. Otherwise, the circles at infinity for universal
coverings of leaves of
the stable and unstable laminations embed in $S_\infty^2$.

\item[e.]
The same results apply to the leaves of any uniform foliation
$\mathcal F$ obtained from $\mathcal F(s)$ by blowing up leaves.
\end{enumerate}
\end{theorem}
\begin{remark} \label{remark: word peano}
This theorem and its proof work equally well if $M$ is
a negatively curved manifold, or simply an irreducible
manifold whose fundamental group is word-hyperbolic
in the sense of Gromov (see for instance \cite{Gromov:wordhyperbolic}.)
However, since a follow-up paper is planned that will prove 
$M$ is a hyperbolic manifold, the statement and its proof are expressed
in terms of hyperbolic 3-manifolds,
for clarity and simplicity.
\end{remark}

\proof
The proof closely parallels the proof in \cite{Cannon:Thurston:peano}.
The basic geometric ingredients are the same as for a hyperbolic
3-manifold that fibers over $S^1$:  uniform spacing between the
leaves and a transverse pseudo-Anosov flow.

Let $l_u$ and $l_s$ denote the 2-dimensional
unstable and stable pseudo-Anosov laminations
($S(I_\nu)$ and $S^*(I_\nu)$ in the notation of section \ref{section: pA}),
and let $l_{uu}$ and $l_{ss}$ denote the one-dimensional
strong unstable and stable laminations, whose leaves are geodesics
on the leaves of $\mathcal F(s)$.
The strategy is first to show that the leaves of $l_u$ and $l_s$
are quasi-isometrically embedded in $\Tilde M$, and that their
circles at infinity converge on $S_\infty^2$. The leaves of these two
laminations will give us enough footholds on $S_\infty^2$ to
pin down the asymptotic behavior of the 
leaves of $\mathcal F(s)$, which bend and wander far more.

We'll use the notation $\Tilde l_u, \Tilde l_{uu}$ \emph{etc.} to
refer to the universal covering laminations in $\Tilde M$. If
$M$ has boundary, then we can think of $M$ as embedded in
the associated complete hyperbolic manifold as a submanifold with
horospherical boundary. We'll denote as $M_+$ the complete hyperbolic
manifold, obtained by gluing horoballs modulo discrete groups to $M$.
Of course, $M_+$ is diffeomorphic with the interior of $M$.
The pseudo-Anosov laminations are contained in $M$ itself; 
the slithering extends to  a slithering $s_+$ of $M_+$, with
associated foliation $\mathcal F(s_+)$. The leaves of $\mathcal F(s_+)$
are complete hyperbolic surfaces; the leaves of
of $\mathcal F(s)$ are obtained by deleting horodisks or their
quotients by $\Z$ (pseudo-spheres).

There is a canonical technique for showing geodesity and quasi-geodesity.
Suppose, for instance, that $\gamma$ is a loop in a Riemannian manifold; how
do you know whether it is the shortest geodesic in its homology class?
The duality between the $L^1$-norm on curves (length) and the $L^\infty$ norm
on $1$-forms gives a necessary and sufficient criterion:
$\gamma$ is minimal in its homology class if and only if there is
a closed $1$-form $\omega$ whose $L^\infty$ norm is $1$ and such that
$|\omega \restrict T_1\gamma| = 1 $.

Similarly, to show that an embedding of path-metric spaces $X \subset Y$ is a
quasi-isometric embedding, a good method 
is to look for a retraction $r: Y \to X$ such that the pull-back by $r$ of the
path-metric of $X$ is a pseudo-path-metric on $y$ that is
quasi-less than the path-metric of $Y$. This translates into a formula
that tries to express the principle:
\begin{multline*}
\exists a > 0 \; \forall x,x' \in X \; \forall y,y' \in Y \\
\left (r(y) = x \; \& \; r(y') = x'\right ) \implies
d(y, y') >  a (d(x,x') - a) .
\end{multline*}
However, the real idea is to keep track of rough distances in $X$
along paths in $Y$,
rather than to analyze all distances at once.   We
can imagine a toll-collector on $X$ who watches the progress of $r$.
Every time $r$ moves further than some threshold $a$ on $X$,
the toll-collector collects a toll, asserting that the
path in $Y$ has gone at least some minimum distance $b$ in
$Y$.  If this assertion were false, the people 
traveling in $Y$ would put in a massive wave of protest; however,
nobody objects when they can travel a long distance in $Y$ without
paying a toll.  The condition on legitimacy of the toll-collection
is logically equivalent to the formula; a
retraction that satisfies this condition is a \df{quasi-isometric retraction}.

\medskip
We'll first analyze the case that $M$ is a closed hyperbolic manifold.  Besides
the hyperbolic metric, we have a second metric that gives a hyperbolic structure
to each of the leaves of $\mathcal F(s)$.  In
the leaf-hyperbolic metric, the projection of a leaf $L$ of
$\Tilde {\mathcal F}(s)$ to any leaf of $\Tilde l_{ss}$ is a quasi-isometric
retraction for the intrinsic geometry of the leaf $L$.  However,  it is not
obvious what happens with this retraction as one varies from leaf to leaf, so
we will construct an alternative.

Instead, we can retract $L$ in $\Tilde M$
to a leaf $g$ of $\Tilde l_{ss}$ by a retraction $r$ that
maps each leaf of $\Tilde l_{uu}$
that intersects $g$ to its intersection point and is monotone
in between, in the sense that $r$ maps the region between two leaves 
of $l_{uu}$ to the interval between their images.

There is an upper bound to
the length of an intersection of $g$ with a gap of $\Tilde l_u$;
clearly, any minimum threshold for assessing the progress of $r$
has to be longer than this minimum length, since
within these intersections, $r$ expands distances
by arbitrarily large factors.  Let
$a$ be a real number greater than the length of
any intersection of any leaf of $l_{ss}$ with any gap of $l_{u}$.
Such a number exists, by compactness.  It follows
that there is a lower bound $b$ to the transverse measure of
a segment on $l_{ss}$ of length $a$, as measured by a pseudo-Anosov
(exponentially shrinking) transverse measure for $l_u$.

The recipe for $r$ on particular leaf of $\Tilde{\mathcal F}(s)$
can be assembled to give a retraction (still called $r$)
of $\Tilde M$ to any leaf $H$ of $\Tilde l_s$, so that each
leaf of $\Tilde {\mathcal F}(s)$ or of $\Tilde l_u$
goes to its intersection with $H$.  Let $p:[0,K] \to \Tilde M$
be any path parametrized by arc length in $\Tilde M$. The
toll-collector on $H$ makes no charge if $|z(p[0,t])| \le 1$
for $t < K$ and if the transverse measure of its projection to
the $\Tilde l_{ss}$ leaf of $r(p(0))$
never exceeds $b$. Otherwise, as soon as
one of these bounds is exceeded, a charge of
\$.25 is imposed, and the accounting is reset. In other words, if
$t$ is the least such time, then 
\[\mathop{\text{charge}}(p\restrict[0,K]) =
\$.25 + \mathop{\text{charge}}(p\restrict[t,K]).
\]

The net toll charged is clearly less than some constant times the arc length of
$p$. Also, the total distance traversed on $H$ by $r(p)$ is clearly less than
some constant times the net toll.  Therefore, $r$ is a quasi-isometric
retraction.

Using the standard principle that in $\Hy^3$ every quasi-geodesic is a bounded
distance from a unique geodesic, it follows that any quasi-isometric
parameterization of $H$ by $\Hy^2$ extends to give an embedding of a closed disk
in $\Hy^3 \cup S_\infty^1$.
A quasi-isometric embedding of $\Hy^2$ in $\Hy^3$ is not usually
within a bounded neighborhood of a hyperbolic plane. Rather, it
can be characterized using quasi-convexity: 
it is equivalent to a uniformly proper image of a topological plane whose
convex hull has bounded thickness---that is, its convex hull separates
$\Hy^3$ into two components whose boundary components are each contained
in a bounded neighborhood of the other.

It is now fairly easy to establish continuity at infinity for a leaf $L$ of
$\Tilde {\mathcal F(s)}$.  Choose a base point $* \in L$. Any geodesic ray $h$
from $*$ on $L$ must have infinite pseudo-Anosov transverse measure for at
least one of the two laminations  $\Tilde l_u$ or $\Tilde l_s$.  This means
that $h$ crosses an unbounded family of leaves, that have eventually empty
intersection with any compact set of $\Tilde M$. This implies that the distance
of their convex hulls from $\*$ tends to infinity, which is the same as saying
that their visual diameter tends to $0$, as seen from $\*$ in the optics of
$\Hy^3$.  In other words, $h$ satisfies the Cauchy condition for convergence in
$\Hy^3 \cup S_\infty^2$.  Furthermore, the regions of $L$ that are cut off by
the leaves of $\Tilde l_{uu}$ and $\Tilde l_{ss}$ that $h$ intersects give a
neighborhood basis for points at infinity of $L$, which shows continuity of the
map of $\Hy^2 \cup S_\infty^1 \to \Hy^3 \cup S_\infty^2$.

\medskip
Now consider the case of a compact 3-manifold $M$ which is a compact
core for a non-compact hyperbolic $3$-manifold $M+$ of finite volume.
We'll first show that the leaves of $\Tilde l_u$  and $\Tilde l_s$
are quasi-isometrically embedded in $\Tilde M$, and then analyze
what this implies.

If $H$ is any leaf of $\Tilde l_s$, we can define a retraction
$r: \Tilde M \to  H$ just as in the closed case, mapping each leaf of
$\Tilde {\mathcal F}(s)$ to its intersection with $H$, mapping
each leaf of $\Tilde l_{uu}$ to its intersection with $H$, and
extending monotonely in between. Notice that this recipe
maps each (horosphere) component of the boundary of the universal cover
to a strip on $H$ between two leaves of $\tilde l_u$.
Since every leaf of $\Tilde{\mathcal F}(s)$ intersects every boundary
component, the image is not bounded above or below.

The laminations $l_s, l_u, l_{ss}, l_{uu}$ are compact,
there is an upper bound to the maximum length
of an intersection of a leaf of $\Tilde l_{ss}$ with a gap of $\Tilde l_u$.
For any real number $a$ greater than this maximum length, there is a lower
bound $b$ to the transverse measure of a segment on a leaf of $l_{ss}$ measured
by a pseudo-Anosov transverse measure for $l_u$.  We can use the same system of
toll-collection as for the compact case. Since $M$ is compact, this
system works, for the same reasons, to show that $r$ is a quasi-isometric
retraction of $\Tilde M$ to $H$, hence that $H$ is quasi-isometrically
embedded in $\Tilde M$.

\medskip
The cusps of $\Tilde M_+$ create logarithmic shortcuts for certain paths in
$M$, and we would hear howls of protest if tried
to extend our system of tolls to $\Tilde M_+$: there is no quasi-isometric
retraction of $\Tilde M_+$ to $H$ that extends $r$. However, we
can use the quasi-geometry of $\Tilde M$ itself, which
is well understood through the work of several people.  

There are several methods that
construct a `sphere at infinity' for a group $G$, with varying
hypotheses; they all
tend to agree in  the simplest situation of a word-hyperbolic group.
For the present circumstance, we can use
a process for compactifying path-metric
spaces of locally bounded geometry, the \df{Floyd compactification},
analyzed by Bill Floyd in \cite{Floyd:thesis}. If $d$ is the 
metric on a space $X$, then we can choose a base point $*$, define a function
$R(x) = d(*,x)$.  For any positive non-increasing
$L^1$ function $f$ on $\R_+$,
there is a metric $d^f$ obtained by measuring path lengths using
the metric $d$ conformally scaled by $f$. If balls of bounded radius
in $d$ are compact, then the metric completion $\Hat X^f$
of the metric $d^f$ is compact.
Let $S^f(X)$ denote the sphere at infinity, $S^f(X) = \Hat X^f \setminus X$.

For some choices of $f$, for instance $f(R) = R^{-2}$, the Lipschitz
class of $d^f$ only depends on the Lipschitz class of $d$. For any such
an $f$, it follows easily that whenever $Q:X \to Y$ is a quasi-isometric
embedding, then $Q$ has a continuous extension
$\Hat Q^f: \Hat X^f \to \Hat Y^f$.

It follows that if $X$ is the universal covering of a compact space
with fundamental group $G$, then $S^f(X)$, up to Lipschitz equivalence,
depends only on $G$. Define $S^f(G)$ to be this sphere.
Floyd showed that when $G$ is the fundamental group of a hyperbolic
$n$-manifold of finite volume, then $S^f(G) = S_\infty^{n-1}$, and
when $G$ is a geometrically finite Kleinian group, then the
limit set of $G$ is the continuous image of $S^f(G)$ under a map which is
usually 1--1, except 2--1 at any rank one cusps of $G$.
\begin{remark}
Note that the price of Lipschitz functoriality of $\Hat X^f$
is infinite Hausdorff dimension, in a case such as for $X$ the universal
cover of a negatively curved surface. The usual metric
for the circle at infinity is obtained as the completion of
$\Hy^2$ using conformal scaling by $f = e^{-R}$, but
the mapping class groups do not act as Lipschitz maps
 in the usual metric.
\end{remark}

This picture is very relevant to our present situation. As a corollary
to the fact that $H \subset \Tilde M$ is quasi-isometrically embedded, 
we obtain a continuous extension 
\[D^2 = \Hat H^f \to \Hat M^f = D^3 , 
\]
where $f(R) = R^{-2}$. This leaves us with the issue of
analyzing any non-injectivity of $H$ at infinity.

Quasigeodesics between points in $\Tilde M$ do not stay within a bounded
distance of each other; this reflects  the fact that
the fundamental group of $M$ is not word-hyperbolic. However, it is
\df{relatively hyperbolic}, relative to its cusp groups. This situation
has been well analyzed, see for example Rich Schwartz's surprisingly
strong classification of the quasi-isometry types of fundamental
groups of cusped hyperbolic manifolds (\cite{Schwartz:quasi-isometry}
and Benson Farb's theory of relatively automatic groups
(\cite{Farb:RelativelyAutomatic}. If we form the quotient $\Bar M$
obtained by collapsing each component of $\boundary \Tilde M$
to a point, then $\Bar M$ is path-hyperbolic, in the sense that
any two quasi-geodesics connecting two points are within
a bounded distance of each other (or the equivalent property, that
a bounded neighborhood of any two sides of a quasi-geodesic
triangle contains the third).  Even though metric
balls of bounded radius are
non-compact in this metric, it still
has the usual sphere at infinity $S^f( \Bar M)$ that is identical with
the usual hyperbolic sphere at infinity
In terms of the geometry of
$\Tilde M$, 
any two quasi-geodesics $g$ and $h$ having the same endpoints,
$g$ is contained in a bounded neighborhood of $h$, 
union any horospherical boundary components that
this bounded neighborhood meets. 
 One can represent 
a quasi-geodesics in $M$ from $x$ to $y$ by a sequence of geodesic
arcs such that any endpoint other than $x$ and $y$ is perpendicular
to a horosphere; then the next arc takes off from some other point
on the horosphere.  In other words, one can think of $\Bar M$ as a
space that turns a quasi-isometrically embedded subset of $\Tilde M$ into
a quasi-convex set.

It follows that any infinite quasi-geodesic whose endpoints at infinity
are identical stays within a bounded distance of some horospherical
boundary component.  

Now we can simply look at quasi-geodesics on $H$ joining its various
points at infinity. In the downward, spreading direction,
we can join two points using two flow-lines of the pseudo-Anosov flow,
connecting them when the distance between them along leaves
of $\Tilde {\mathcal F}(s)$ decreases $1$.  The pseudo-Anosov
 transverse measure for $\Tilde l_u$ between these flow-lines
grows to infinity, which implies that these ends are not
contained in a bounded neighborhood of any single cusp.

There is a unique point at infinity in the upward, contracting direction.
If $H$ is the face of a cusp gap of $\Tilde l_s$, then a closed loop
on $H$ is homotopic to a cusp; its two endpoints are identified, and
all other endpoints necessarily are distinct.
In other words, $H$ makes a figure \textsf{8} on $S_\infty^2$.
We can surger this cylindrical  `accidental parabolic' leaf
$H/\Z$ into two `deliberately parabolic' pseudospheres $H_1$
and $H_2$, using a saw-blade from
of section \ref{section: pA}) for the surgery. The universal covers of the
resulting pseudospherical leaves have completions that are disks.

In the generic case when $H$ is not the face of some cusp gap of $\Tilde l_s$,
then in the upward direction,
its $\Tilde l_s$ transverse measure to any cusp is non-zero on
any leaf $\Tilde {\mathcal F}(s)$, so it tends to infinity as $z \to \infty$.
Thus it does not remain in any bounded neighborhood of any cusp,
so it is not identified with any point at infinity in the downward
direction. In this case, $\Hat H^f \to S_\infty^2$ is injective.

\medskip
Now consider a leaf $L$ for $\Tilde {\mathcal F}(s)$. Let $*$ be
a base point on $L$, and consider any geodesic ray $h$
emanating from $*$ in the hyperbolic metric of $L$.

If $h$ does not
tend to a cusp of $L$, then it crosses infinite transverse pseudo-Anosov
measure for at least one of the two laminations $\Tilde l_{uu}$ 
or $\Tilde l_{ss}$.  This
implies that the ray enters (and stays in) 
half-spaces cut off by leaves of $\Tilde l_u$ or $\Tilde l_s$ that are
infinitely far from $*$ in $\Bar M$. Since these neighborhoods are
quasi-convex in $\Bar M$ and arbitrarily distance,
they are arbitrarily small in the completion; it follows
that $L$ converges at infinity and is continuous at the endpoint of
such an $h$.

If $h$ tends to a cusp, it only traverses a finite total
transverse measure for either $\Tilde l_u$ or $\Tilde l_s$. 
We need a slightly different construction: we can use
the doubly-infinite sequence of leaves of $\Tilde l_{u}$ and
$\Tilde l_s$. Small neighborhoods for the endpoint of $h$ in the
domain $D^2 = \Hat{ \Tilde L}^f$ can be cut out by using one half-leaf on each
side of $h$  of either
$\Tilde l_{uu}$ or $\Tilde l_{ss}$,
together with a portion of a horocycle. In $\Hy^3 = \Tilde M_+$, these
map into neighborhoods cut off by two of the deliberately
parabolic half-leaves from $\Tilde l_u$ or $\Tilde l_s$ that
were formed by surgeries, together with a strip that joins them on a
horosphere.  To see that sets of this form
shrink in size to the parabolic point, we can use the fact that
the system of all the deliberately-parabolic half-leaves is invariant
by the $\Z+\Z$ that stabilizes the cusp. In other words, they come from
a finite set of parallel pseudospheres in $\Hy^3 / \Z+\Z$. It follows
that all but a finite set of their limit circles in $S_\infty^2$
have size less than a given constant $\epsilon$. The assemblage
of two surger ed half-leaves plus a strip on the horosphere limits
at infinity to a figure \textsf{8} formed by combining small wings from
two limit figure \textsf{8}'s (a non-topologist would be more
likely to call them sausages with ends joined)
of leaves of $\Tilde l_u$ or $\Tilde l_s$. These
wings becoming arbitrarily small. This shows convergence of
$D^2 = \Hat L^f$ near a parabolic point on its boundary.
\endproof

\section{Anosov flows and extended convergence groups}
\label{section: Anosov}

Much inspiration for the present study came from
S\`ergio Fenley's interesting analysis 
of Anosov flows on $3$-manifolds (\cite{Fenley:skew}.) Fenley
developed a surprisingly strong theory for certain Anosov flows and their
associated foliations.   From Fenley's results, interpreted in terms of
slitherings, a beautiful and suggestive picture emerges, a picture
that suggests there is much more that is yet to be understood.

An \df{$\R$-covered} foliation is a foliation such that the 
space of leaves in the universal cover is homeomorphic to $\R$.
An Anosov flow $\psi_t$ on a $3$-manifold $M$ is called
$\R$-covered if its stable foliation $\mathcal F_s$ is $\R$-covered.
Fenley proved that an Anosov flow of a $3$-manifold is $\R$-covered,
then it has one of two types.  The first type is the \df{product} type,
when every leaf of $\Tilde{mathcal F}_u$ intersects every leaf
of $\Tilde {\mathcal F}_s$; this
happens if and only if $\psi$ is the suspension of an Anosov diffeomorphism
of $T^2$. 

The second type is that of a \df{skew $\R$-covered} Anosov flow. In this case,
$\tilde M$ can be mapped surjectively
to the diagonal strip $|x-y|<1$ in the plane
so that the preimage of any point is a flow-line of $\Tilde \psi$,
the preimage of any horizontal lines is a leaf of
$\Tilde {\mathcal F}_u$, and the preimage of any vertical
line is a leaf of $\Tilde {\mathcal F}_s$. 

The geodesic flow for a hyperbolic surface, illustrated
in figure \ref{figure: circle at infinity}, is the primordial example.
In that figure, the surface of the cylinder is divided,
like a mailing tube, into two bands that wrap diagonally around it.
The flow-lines of $\Tilde \psi$ are
the horizontal lines inside the cylinder; if stable leaves are projected
in one direction, they map to the foliation of one of the
bands by vertical lines, while projecting in the other direction maps
each stable leaf to the diagonal foliation of the other band, where
the leaves wrap with slope $1/2$. (This is the same $1/2$ as in `an angle
inscribed in a circle is $1/2$ the central angle').   The unstable
foliation is obtained by rotating the picture $180\degree$ about
its vertical axis; this gives two foliations of each strip. Either
strip gives a good initial model of a skew $\R$-covered
Anosov flow. 

Many further examples of skew $\R$-covered
Anosov flows that can be constructed by Dehn surgery along closed
trajectories of Anosov flows. On the boundary of a regular neighborhood
of a closed trajectory, there are distinguished
closed curves, coming from the intersection with the stable and unstable
leaves of the trajectory.  It has been known for some time that
any surgery obtained by re-attaching the regular neighborhood by
a diffeomorphism that preserves these curves (that is, 
by a power of the Dehn twist about one of these curves)
yields another 3-manifold with an Anosov flow.  
Fenley showed that if the original flow is the suspension of an Anosov
diffeomorphism of $T^2$, then any surgery that uses consistently-oriented
Dehn twists along any collection of closed orbits yields a skew $\R$-covered
example.  
This fits with a construction 
of Hedlund and Morse (?) in which they constructed sections for the
geodesic flow for any hyperbolic surface in the complement of certain
systems of closed geodesics (in $T_1(M^2)$
homeomorphic to a multi-punctured torus---in other words, the geodesic
flows are obtained by Dehn surgery from suspensions of Anosov diffeomorphisms.
This construction also shows that the geodesic flow for any oriented
hyperbolic orbifold is obtained by surgery from an Anosov suspension.

The simplest case where everything is orientable
is the figure eight knot complement, which fibers over the circle
with fiber a punctured torus glued by
$\bigl [
\begin{smallmatrix}
2&1\\1&1
\end{smallmatrix}
\bigr ]$.
In \cite{Thurston:GT3M} it was shown that several of the Dehn fillings
give Seifert fiber spaces. Every manifold obtained in this entire row of
Dehn fillings of the figure knot has a skew $\R$-covered Anosov flow,
according to Fenley's analysis; the cases that are Seifert fibered cases are
examples of Hedlund's construction. 

Furthermore, Fenley showed that for any skew $\R$-covered Anosov
flow, there is an orientation defined by the structure so that any
positive foliation-consistent surgery along orbits yields another
skew $\R$-covered example.

\medskip
A key theme of \cite{Fenley:skew} is that every
automorphism of the double foliation of the diagonal strip is
periodic.  For any stable leaf, there is a lowest unstable leaf that
doesn't intersect it; for every unstable leaf, there is a lowest
stable leaf that doesn't intersect it. This gives a canonically-defined
equivalence relation on the set of stable leaves that is necessarily
preserved by any automorphism. The quotient of the equivalence relation
is a circle.  In other words,
$\mathcal F_s$ is the foliation of a slithering
of $M$ around $S^1$---we can use the fact 
that $\Tilde {\mathcal F}_s$ is equivalent to a product foliation
$\R^2 \times \R$ (\ref{proposition: R-covered R-slitherings}) to 
check that the map $\Tilde M \to S^1$ is a fibration.

One can picture this as in \ref{figure: circle at infinity}: 
glue together two copies of the skew strip,
with stable and unstable directions interchanged, to form
a cylinder. The equivalence classes correspond to the circle's worth
of vertical lines---the intersection of one vertical line with
one of the strips describes a set of equivalent stable leaves.
The gluing map, which interchanges the two sides of the strip by
going either straight up or straight to the right, we may as
well call $\sqrt Z$, since $\sqrt Z \circ \sqrt Z $ represents
$Z$.  Here's what the picture says so far when translated further into 
the logic medium:

\begin{proposition} \label{proposition: skew identification 1}
Let $\psi$ be a skew $\R$-covered Anosov flow on a 3-manifold $M$.
Then there are two slitherings $S_s$ and $S_u$ of $M$ around $S^1$,
where $\mathcal F(S_s) = \mathcal F_s$, and $\mathcal F(S_u) = \mathcal F_u$.

The circle at infinity bundle for the leaves of  $\mathcal F_s$ is
isomorphic to the slithering circle bundle for $\mathcal F_u$.
and vice-versa. The isomorphism is obtained by 
gluing two copies of the closure of the skew strip using $\sqrt Z$, and
attaching it equivariantly to 
$\Tilde M$  to form a solid cylinder.
When the `short' leaves of one of the strips are collapsed,
the `long' leaves of the second
strip join to become the circles at infinity for
one of the foliations, while the `short' leaves of the second strip
join to become lines representing the quasi-isometric identification 
of the leaves at infinity.
\end{proposition}
\proof
Since  the strip parameterizes flow lines of $\Tilde \psi$, and
each of these flow lines is $\R$, we can sew two copies of the open strip
to $\Tilde M$, each point attached to one of the two ends of its flow line,
and we can adjoin lines to serve as common edges for the two strips.
A leaf of an Anosov flow automatically looks like the hyperbolic plane
foliated by geodesics emanating from one point, and the collapsing
maps of the strips just describe this geometry.

The only actual issue is the comparison between
the identification of circles at infinity according to bounded distances
between two geodesic rays and the identification defined using the
skew-Anosov structure.  Given two nearby leaves $L_1$ and $L_2$
of $\Tilde {\mathcal F}_s$,  the matching is forced along
most of their circles at infinity, wherever a leaf of $\Tilde {\mathcal F}_u$
meets both $L_1$ and $L_2$. Consider the cylinder formed by collapsing
the short leaves of one of the strips, which is like the core of 
a roll of paper towels.  We know that the horizontal
circles are circles at infinity for the leaves, and we know that
the vertical foliation gives the correct identification of the circles
everywhere except possibly on the boundary of the strip $/Z$. But it
is clear (and easy to prove) that 
the only extension of the vertical-line foliation across the missing line
is the foliation by vertical lines.
\endproof

The four actions (slitherings and leaves at infinity of the two foliations),
which the proposition says are actually only two different actions,
are in fact isomorphic
to each other, since $\sqrt Z$ conjugates one action to the other.
Another way to rephrase the picture so far is in terms of the
space $P$ of ordered pairs of distinct points on a circle.
$P$ is homeomorphic to an annulus, and $\tilde P$ is the diagonal
strip, with two foliations coming from the two coordinates. An
element of $P$ represents a geodesic in the hyperbolic plane, and
we can think of $\pi_1(M)$ as acting on the space of geodesics.

\subsection{Extended convergence groups}
\label{subsection: extended convergence}
Now it's time for a third foliation $\mathcal F_p$ to enter the picture.
The points on a geodesic $\overline{us}$ in $\Hy^2$ can be parametrized 
in terms of a third element $p \in S_\infty^1$ that
counterclockwise from $s$: the foot of the perpendicular
from $p$ to $\overline{us}$ gives a $1-1$ correspondence.

Let $T$ be the space of counterclockwise ordered triples of distinct
points on $S^1$. We can complete $T$ by admitting degenerate cases
when some or all the elements of the triple coincide; however, we
remember their limiting order, so that when all three
clumping in the order $usp$ we distinguish it
from clumping in the order $spu$, even when the points clump
at the same place.
When we do this, we obtain a solid torus $\Bar T$ whose boundary is
divided by the three $(1,1)$ curves into three annuli.

\begin{definition} \label{definition: extended convergence}
A \df{convergence group} is a subgroup of $\Homeo_+ S^1$ that
acts properly discontinuously on $T$, where $\Homeo_+$ denotes
the orientation-preserving subgroup.

An \df{extended convergence group} is a subgroup
of the universal covering group $\widetilde{\Homeo_+ S^1}$ (consisting
of periodic homeomorphisms of $\R$) that acts properly discontinuously on
$\Tilde{\Bar T}$.
\end{definition}
We can coordinatize $\Tilde T$, when convenient, as the set of ordered
triples of real numbers $(u,s,p) $ where $u < s < p < u + 2 \pi$.

\begin{definition} \label{definition: total foliation}
A \df{total foliation} for an $n$-manifold is a collection of
$n$ codimension one foliations, locally equivalent to the $n$ foliations
of $\R^n$ parallel to the coordinate axes.
\end{definition}

Detlef Hardorp (\cite{Hardorp:total}) proved that every 3-manifold
admits a total foliation. His construction makes free use
of Reeb components; as far as I know, there has been little 
investigation of \df{taut total foliations}, that is, total
foliations such that the three codimension one foliations are taut.

The annulus $D$ and the solid torus $T$ come equipped with total
foliations. The group of automorphisms of the total foliation of $D$ or of
$T$ is isomorphic to the group of homeomorphisms of $S^1$. Similarly,
the groups of automorphisms of the total foliations of $\Tilde D$
and of $\Tilde T$ are isomorphic to $\widetilde \Homeo(S^1)$.
For any  extended convergence group $\Gamma$, the three-manifold
$\Tilde T / \Gamma$ comes with a built-in taut total foliation.

Here is a collection of basic properties of extended convergence groups:
\begin{proposition} \label{proposition: basic extended convergence}
Let $\Gamma$ be an extended convergence group, and let $M = \Tilde T / \Gamma$
be its quotient three-manifold. 
\begin{enumerate}
\item[i.]
$M$ has three slitherings $S_u$, $S_s$ and $S_p$ around $S^1$ whose
foliations form the taut total foliation of $M$. The cartesian
product $S_u \times S_s \times S_p: \Tilde M \to S^1 \times S^1 \times S^1$
is the covering map $\Tilde M \to  T \subset T^3$.
\item[ii.]
There is a homeomorphism $Z^{1/3}: \Tilde M \to \Tilde M$ that
commutes with deck transformations of $\Tilde M \to M$,
cyclically permutes the slitherings and whose cube acts as $Z$
on the space of leaves of all three slitherings.
\item[iii.]
The leaves of the three total foliations are of hyperbolic type,
that is, there is a Riemannian metric on $M$ that restricts to
a hyperbolic metric on the leaves of one of the slitherings.
\item[iv.]
An element of $\pi_1(M)$ that is space-like for one of the slitherings
is space-like for all three.
The action of any space-like element of $\pi_1(M)$  on $S^1$ (the circle
where the three points lie) either has one fixed point (parabolic case),
or has two fixed points where one is attracting and one repelling
(hyperbolic case).
If $M$ is compact, then the parabolic case does not occur.
\item[v.] If $M$ is compact, or if it is
the interior of a compact manifold to which the slitherings extend, then
the circle at infinity for the leaves of any of the
three foliations can be identified with the original $S^1$
(of the triples) so that matches the actions of $\pi_1(M)$.
\item[vi.] The completion of $\tilde M$ by the circles at infinity
for the leaves of $\Tilde{\mathcal F}(S_s)$ is homeomorphic to the
solid cylinder obtained from the triangular prism
$\Tilde{ \Bar T}$ by collapsing its two faces $u=s$ and $s=p$
along the lines where $s$ is constant.
\end{enumerate}
\proof
Part (i.) is really a rephrasing of the definition of an extended 
convergence group.

Part (ii.) comes from the symmetry of the three elements of
the triple on $S^1$. In terms of the coordinates $(u,s,p) \in \R^3$
where $u < s < p < u + 2 \pi$, the map is $(u,s,p) \to (s,p,u+2\pi)$.

Part (iii.) is almost vacuous when $M$ is non-compact unless 
further conditions were put on the metric: for a foliated 3-manifold
where no leaf is contained in a compact subset,
it's easy to modify the metric near infinity to
make all leaves of hyperbolic type. When $M$ is compact, the cocompactness
of the action of $\Gamma$ on $\Tilde T$ implies that any pair of points
$u < s < u+2\pi$ in $\R$ can be
that squeezed arbitrarily close together by $\pi_1(M)$.
No measure on $\R$ can possibly be invariant under this action.
From Candel's theory of uniformization of surface laminations
(\cite{Candel:Uniformization}), it follows that the leaves of
the foliations are conformally hyperbolic, for otherwise there would
be an invariant measure.

For part (iv.), the homeomorphism $Z^{1/3}$ that commutes with the
action of $\pi_1(M)$ shows that space-likeness is equivalent for the
three foliations.

The action of a space-like element of 
$\widetilde{\Homeo S^1}$ on $\R$ is to fix all points that cover
its fixed points on $S^1$; a space-like element that fixes a
three or more points on $S^1$ therefore fixes
some element of $\Tilde T$, and cannot be part of a properly discontinuous
group. (Note that all elements of $\widetilde{\Homeo S^1}$ have infinite
order. In particular, the elements of
$\widetilde {\Homeo S^1}$
that cover torsion elements of $\Homeo S^1$ are time-like and
of infinite order).

A little further thought shows that if $\gamma \in \Gamma$
is a space-like element with two fixed points $a$ and $b$, then iterates
of $\gamma$ take all points on the circle except for one of the points,
say $b$, toward the other fixed point ($a$)
(and $\gamma^{-1}$ does the reverse). Otherwise, $\gamma$ would
map the two intervals between $a$ and $b$ in the same sense,
say counterclockwise, and the action of $\gamma$ on $\R$ would be similar,
with $\gamma(x) \ge x$.
In this situation, one can find a compact arc in $\Tilde T$ whose images
under iterates of $\gamma$
get hung up in a compact set, and never escape to the boundary: if
you think of $\Tilde T$ as triples of real numbers contained in
an interval of size less than $2 \pi$, then keep the two outside
points above and below one of the fixed points of $\gamma$, while
letting the third point cross the fixed point. This shows that
every space-like element is either of hyperbolic or of parabolic type.

When $\Gamma$ is cocompact, we can
rule out parabolic elements more reasoning of a similar nature.
Indeed, suppose there were a parabolic element $\gamma$ of
a cocompact extended convergence group $\Gamma$.
Let $a \in \R$ be a fixed point,
let $a < b < a+2\pi$ be another point, and let $c = \gamma(b)$.
Cocompactness would imply 
that some element of $\Gamma$ would take
$b$ and $c$ close together in $(a, a+2\pi)$ while keeping them
keeping them far from $a$. Let $\gamma_i$ be a sequence of
conjugates of $\gamma$
by such transformations. The sequence $\gamma_i$ is necessarily
unbounded in $\Gamma$; proper discontinuity therefore requires
that for 
any two compact sets $K$ and $L$  in $\Tilde T$, $\gamma_i(K)$
is eventually disjoint from $L$. But the sequence $\gamma_i$
has qualitative behavior very similar to a homeomorphism 
with two fixed points but the wrong dynamics---an arc of triples can
be constructed  that gets hung up in a compact subset of $\Tilde T$.

To establish part (v.), assume first that $M$ is closed, and
let $L$ be a leaf of $\Tilde{\mathcal F}(S_s)$.
By proposition
\ref{proposition: space-like generators}) $\Gamma$ contains non-trivial
space-like elements. Thus,
we can compare the two circles $S^1$ (containing the original
triples) and $S_\infty^1(L)$ by looking at fixed points of space-like
elements of $\Gamma$, using the fact that all leaves are 
$\Tilde {\mathcal F}(S_s)$ are quasi-isometrically equivalent.
If $\gamma$ is a space-like element, then it fixes some leaf $L'$ in
the slab between $L$ and $Z(L)$, and it necessarily acts as
a hyperbolic element
on this leaf (using the fact there is a lower bound to injectivity
radius of leaves.) Any quasi-geodesic invariant by $\gamma$ converges in
one direction to the attracting fixed point of $\gamma$,
in the other direction to the repelling. We can conjugate any space-like
$\gamma$ by any element of $\pi_1(M)$; it is an easy exercise to see
that the fixed points of conjugates must be dense in $S_\infty^1(L)$.

The attracting fixed points of space-like elements
on  $S^1$ and $S_\infty^1(L)$ inherit circular orderings which
can be reconstructed from the topology of $M$ by looking at orientation
information coming from intersections of closed
geodesics with annuli in $M$
illustrated in figure \ref{figure: Lambda} and used (without
orientations) for the construction of the linking series $\Lambda$.
This gives a 1-1 circular-order-preserving
identification of a dense set of points on the two circles, which
therefore extends to a homeomorphic identification.

When $M$ is the interior of a compact manifold to which the slitherings
extend, as usual we define the quasi-isometry type of leaves by gluing horoball
quotients to the torus boundaries. The easiest argument in this
case is to use the cusps on leaves to define the identification
of the two circles. Cusps are indelibly printed on
each leaf as components of its intersection with a neighborhood of infinity,
so the circular ordering has an immediate topological (not just
quasi-isometric) definition.

Part (vi.) is a reworking of part (v.) to picture it in terms of $\Bar T$
and $\widetilde {Bar T}$.
Given a space-like element $\gamma \in \Gamma$, 
the effect of applying iterates of $\gamma$ to an element of 
$T$ is to take at least two of the three elements of the triple
to the attracting fixed point of $\gamma$ on $S^1$.  Therefore,
for any quasi-geodesic in $\Tilde T$ that is contained between
two leaves $L$ and $Z(L)$ and is
invariant by $\gamma$
has limit set contained in
the subset of $\Bar T$ where at least two of the three
coordinates have value equal to one of the points in $\R$ that covers
the convergent fixed point of $\gamma$ on $S^1$.  Limit sets for
space-like elements are dense; possible limits of other quasi-geodesics
are similar, and can be deduced either from how they are sandwiched between
group-invariant quasi-geodesics. (There is another approach for (v.) and (vi.)
that perhaps helps clarify the picture. One can look
directly at space-like quasi-geodesics in a $z$-bounded slab of $\Tilde M$,
and use the sequence of fundamental domains that they intersect. This
sequence is labeled by a sequence of elements of $\Gamma$, which
a space-like quasi-geodesic for the group.  It is not hard to
see that the behavior of this sequence is similar to the special
case of iterates of a single space-like element: for any particular
large $i$, all of $S^1$ goes near a particular `attracting'
point, all of $S^1$ except a short interval $J_i$; however, the intervals
$J_i$ are not located in any consistent place.)

When we restrict to a leaf where $s$ is constant, the two faces
$u=s$ and $s=p$ therefore represent a single point at infinity for $L$,
while the remaining face of the prism sweeps transversely across
the possible limit sets of closed space-like quasi-geodesics, filling out
the rest of the circle.
\endproof
\end{proposition}

\begin{corollary} \label{corollary: skew triples}
The quotient of $\Tilde T$ by any orientation-preserving
cocompact extended convergence group is a $3$-manifold with
a total foliation such that a vector field tangent to
the intersection of any two of the three
codimension one foliations is a skew $\R$-covered Anosov flow.

Conversely, every skew $\R$-covered Anosov flow has this form.
If $s$ is a slithering of a closed orientable
$3$-manifold $M$ around $S^1$, then $\mathcal F(s)$ is the
stable foliation of an Anosov flow
if and only if the associated representation
of $\pi_1(M)$ in $\widetilde \Homeo(S^1)$ is a cocompact extended convergence
group. In that case, $M = \Tilde T / \pi_1(M)$. 
\end{corollary}
\proof
The logic is easy, based on what we know. Given a skew $\R$-covered
Anosov flow $\psi$ on a $3$-manifold $M$, map $\Tilde M$ to $\Tilde T$,
the map to $\Tilde D$ to to give two of the three coordinates $u(x)$ and $s(x)$
for $x \in \Tilde M$.
For the third coordinate, we choose hyperbolic structures for the leaves
of $\Tilde {\mathcal F}_s$, and let
$p(x)$ be the leftward endpoint on $S_\infty^1$ of
the perpendicular to the $\Tilde \psi$-flow line through $x$,
in the geometry of the $\tilde {\mathcal F}_s$-leaf of $x$.
The orientation information comes because of the direction of skewing:
in other words, the strip $\Tilde D$ has an orientation, and
the orientation of the flow gives an orientation for $M$.
On any particular flow line of $\Tilde \psi$, there is an interval's
worth of choices for $p$, which is the right information to determine
a point on the fiber of the map $\Tilde T \to \Tilde D$ that forgets $p$.  
Note that this coordinatization is not smooth; however,
smoothness is not a critical issue when we have the strong structure
provided by a total foliation.  The action of $\pi_1(M)$ on $\Tilde T$ is
properly discontinuous and cocompact, because we have produced
an isomorphism to the action of $\pi_1(M)$ by deck transformations,
so $\pi_1(M)$ is an extended convergence group.

If $s$ is a slithering of $M$ around $S^1$ and if $\mathcal F(s)$
is an Anosov foliation, it is an $\R$-covered Anosov foliation,
so it follows from \cite{Fenley:skew} that it is skew $\R$-covered,
and therefore it acts on its space of leaves as an extended convergence
group.

What remains is to establish that the quotient of a cocompact extended
convergence group $\Gamma$ has Anosov flows as stated.
Given our understanding of the circles at infinity of leaves, we
can see this just by looking. If we look at two of the foliations
of $\Tilde T$, say the foliations where $s$ is constant and where
$u$ is constant, then the leaves of the
1-dimensional foliation $s = u =  \text{constant}$ restricted to
a leaf $s = \text {constant}$ converge to a single point in
one direction (because of collapsing of two faces of the prism), and
they converge in the opposite direction on a leaf \(u = \text{constant}\).
It is easy to see that convergence has to be exponentially fast
by looking at a hyperbolic metric for the leaves of one of the
foliation, and applying general considerations of compactness: therefore,
a vector field tangent to this foliation is an Anosov flow.
\endproof
\begin{remark}
Notice that if we were to collapse all three faces of the prism
$\widetilde{\Bar T}$ along their `short' directions, the boundary
would collapse to a circle and $\Bar T$ would collapse to an
uncoated lens.  However, if $\Gamma$ is not a convergence group, then
the quasi-isometric distance between leaves at
different levels of $\Tilde T$ goes to infinity with height, so just
like the original prism, this
lens must be interpreted used cautiously for understanding the
quasi-geometry of $\Tilde T$.
If the foliations are transversely pseudo-Anosov,
a transverse pseudo-Anosov flow $\phi$ gives a connection that can
be used to define a compactification of $\Bar T$, by re-mapping
the interior of the prism to a compact triangular prism, 
fixing one leaf $s = \text{constant}$, then mapping
flow lines of $\phi$ to parallel lines that terminate at parallel
end faces. Now when we collapse the three short directions of the
rectangular sides of the prism, we get a coated lens (with its
faces).  As quasi-geodesics,
the stable leaves of the transverse pseudo-Anosov flow all
collapse at the top of the lens, and the unstable leaves all converge
at the bottom.    Collapsing the leaves of these foliations
on the top and bottom of the lens yields a ball; the rim of the lens
becomes the sphere-filling curve of section \ref{section: Peano curves}.
\end{remark}

\bigskip
Here is another representation for these same manifolds, adding
a twist of contact geometry, whose primordial example is a cotangent
sphere bundle (or unit cotangent bundle, if one prefers to use a metric).
The steady rotation of flow-lines of the Anosov flow, as one goes
transversely to the leaves of $\mathcal F_s$, suggest a contact structure.
This can be used to compile several structures all into one picture:

\begin{figure}
\centering
\begin{minipage}{.45\textwidth}
\includegraphics[width=\textwidth]{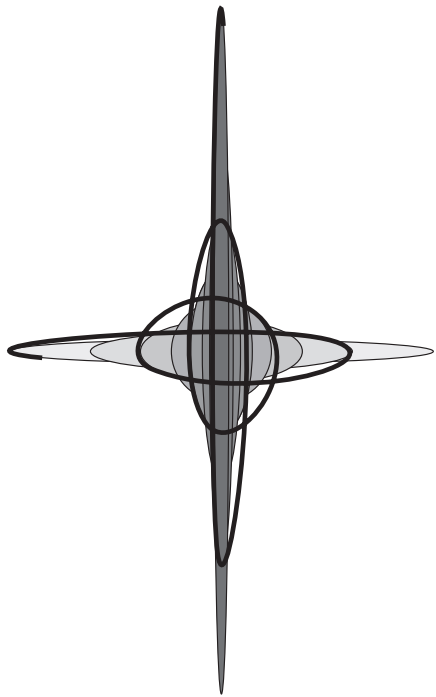}
\end{minipage}
\begin{minipage}{.55\textwidth}
\caption[Anosov foliations come from scribblings]{
\label{figure: scribbles}
The universal cover of a manifold with a skew $\R$-covered Anosov foliation 
has an immersion in $\widetilde{ T^*(\Hy^2)}$ that is invariant
by the derivative action of a cocompact properly discontinuous
subgroup of diffeomorphisms of the hyperbolic plane. In any fiber,
the immersion would look something like this. There is a 1-parameter
family of coaxial ellipses that describe the conformal structure of
each leaf. The immersion is tangent to the unit ellipse for its hyperbolic
metric.}
\end{minipage}
\end{figure}

\begin{theorem} \label{theorem: scribbles}
The action of an extended convergence group on $\Tilde T$
is conjugate to the action of a subgroup of $\widetilde {\Diff}_+(\Hy^2)$
of uniformly Lipschitz diffeomorphisms on $\widetilde {\TS}^*(\Hy^2)$,
where the action is the universal cover of the derivative.
Conversely, any subgroup of the
universal covering of the group of uniformly
Lipschitz diffeomorphisms of $\Hy^2$ which acts properly discontinuously
on $\widetilde{TS}^*(\Hy^2)$ is conjugate to an extended convergence
group acting on $\Tilde T$.

Furthermore,
if $M$ is a closed $3$-manifold with a skew $\R$-covered Anosov foliation,
then $\Tilde M$ can be immersed in $T^*(\Hy^2)$ (probably not
very smoothly), in a way that is invariant by the derivative
of a representation of $\pi_1(M)$ in $\widetilde{\Homeo}(\Hy^2)$.
\end{theorem}
\proof
Theorems \ref{theorem: convergence groups} and
\ref{theorem: transverse pseudo-Anosov} can be interpreted in
the context of extended convergence groups, in the following form:
Given an extended convergence group, we choose one of the
three (isomorphic) slitherings, and apply theorem
\ref{theorem: transverse pseudo-Anosov} or \ref{theorem: convergence groups},
to obtain a flow $\phi$ transverse to the foliation of the slithering.
For present purposes, we just need a flow that serves as a connection
for the slithering, so if necessary, we perturb $\phi$ to a
flow $\phi'$ that is a smooth and remains a connection.

We can now represent $Z$ by the homeomorphism having the
required action on leaves of $\Tilde M$, and isotopic to the identity
along flow lines of $\phi'$.   This homeomorphism is likely not
to be smooth, but it acts smoothly on the space of flow-lines.
This gives a representation of $\pi_1(M)$ as a group of homeomorphisms
of the universal cover of any leaf.  We can lift this representation
to the universal covering group of the group of homeomorphisms
of the leaf, using the homotopy information that comes from
the slithering, since the action on cotangent circles is the same
up to homotopy with the action on the circle at infinity, which action
is equipped with a lifting to an action on $\R$.
\endproof

\section{Preview and questions}

The circle at infinity for the leaves of a slithering is a 
particularly well-behaved instance of a general construction for
a universal circle-at-infinity for the leaves of any taut foliation
of a 3-manifold. In general, the universal circle can be thought
of as defined by a foliation transverse to the fibers of $TS(\mathcal F)$.
It is not homeomorphic to the circles at infinity defined by the
geometry of individual leaves. Instead, it is a collation of
the circles for all leaves into one master circle:  the circle
at infinity for any particular leaf is obtained as a monotone (but
not strictly monotone) image of the master circle. These universal
circles will be constructed and analyzed in \cite{Thurston:circlesII}
and they will be used to construct genuine essential laminations
transverse to the leaves of any taut foliation of an atoroidal
3-manifold. 

Harmonic measures for foliations constructed by Garnett (\cite{MR84j:58099})
are very helpful in understanding the geometry of leaves of taut foliations.
They can be used to show that on any leaf $L$,
in `most' directions at infinity (in some sense, where the exact meaning of 
`most' depends on $L$)
the holonomy keeps a definite packet of nearby leaves within a bounded
distance, and under many circumstances, makes them converge toward $L$.
Anosov foliations are a particularly clear instance of this, where
on any leaf, in all directions at infinity except one, the flow-lines 
diverge, but nearby leaves converge. The general picture is similar
to this, except that there is often a dense set of exceptional
directions where nearby leaves diverge.  The exceptional set of
directions has measure $0$ in any foliation such that every leaf is dense.

There has been a long history of a need for a widely-applicable geometric
theory of \df{universal Teichm\"uller space}, that is, the space
of hyperbolic structures on $D^2$ \textit{rel} boundary, subject
to some constraint on the geometry. This has been a key issue
in studying iterated rational maps of the Riemann sphere, and it
is also a key issue in the topology of three-manifolds.   Of course,
there are also many interesting unresolved issues concerning the
geometry of ordinary finite-dimensional Teichm\"uller spaces.

I believe that three-manifolds that slither around $S^1$ provide a nice
attainable testing ground, for refining some of our understanding
about hyperbolic geometry and Teichm\"uller geometry.  In general,
given a compact space with a 2-dimensional lamination that has
a hyperbolic leaves, one can study the Teichm\"uller
space for the leaves---what are the possible hyperbolic metrics, up
to isometry? 

Associated with any lamination with 2-dimensional
hyperbolic leaves, there is an associated 3-dimensional lamination
with 3-dimensional hyperbolic leaves, modeled on
$\Hy^3 \supset \Hy^2$.  The 3-dimensional laminations have interesting
deformation spaces of their own; these are generalization of 
quasi-Fuchsian groups.

In the case of a lamination $\lambda$ embedded
in a 3-manifold $M^3$, one can go further, and incorporate hyperbolic structures
on the gaps into a foliation of 4-manifold $N^4$ with 3-dimensional hyperbolic
leaves; $M^3$ is embedded in $N^4$ as a spine, in a way that leaves
of the foliation of $N^4$ intersect $M^3$ as the leaves and gaps. (One
should first blow up any isolated leaf of $\lambda$ to a band of parallel
leaves, to make sure that $N^4$ will be Hausdorff.)

There are fairly natural ways to define a relaxation process on
the deformation space of the 3-dimensional hyperbolic foliation,
to try to bring nearby leaves isometrically closer. Actually,
similar processes can be defined on the leaves of a foliated
3-manifold;  some of these processes conjecturally
should tend to a limit that is analogous to a geodesic in one of the metrics
for Teichm\"uller space, yielding a transverse pseudo-Anosov flow
under fairly general circumstances.  (This would generalize
`curve shortening' as carried out by Bers in the special case of a
surface fibration  over $S^1$, \cite{Bers:curve-shortening}.) 

I think it is likely that a relaxation process can be defined
for the three-dimensional hyperbolic foliations that converges to
give a geometric decomposition for $M$, usually by converging
to a foliation where all leaves are actually isometric.  What appears
to happen is that as that an appropriate relaxation process
makes the leaves `rotate' in $\Hy^3 \times $ Teichm\"uller space, so that
the up and down $\Hy^3$-directions turn toward neighboring leaves above
and below, ultimately converging to be isometric if the lamination
is irreducible in the sense of not admitting transverse essential tori.

The case of a lamination
consisting of a finite number of incompressible surfaces is
simply a translation of a Haken manifold 
into this language, and the proofs for Haken manifolds
show that the relaxation process converges  in this case.

The next case will be 3-manifolds that slither around $S^1$. I am
planning a paper to prove geometric 
convergence of a relaxation process defined in the pseudo-Anosov
case by $(Z,Z^{-1})$.  A key ingredient is that 
as one iterates, the quasi-isometric distance between leaves never
increases, just as the `skinning maps' decrease the Teichm\"uller metric
in the Haken situation.  In fact, because of the uniformness of the
quasi-geometry of leaves in this situation, Curt McMullen's proof of the
Theta conjecture shows that distances between leaves actually contracts.
Geometric estimates similar to previous cases will imply geometric
convergence, yielding a hyperbolic structure for
$M^3$. This scheme is a generalization of
the construction for hyperbolic structures
on mapping tori of pseudo-Anosov homeomorphisms.

I think it likely that further study will eventually
show convergence (of some version of this process)
in generality, for 3-manifolds with taut foliations
or essential laminations.

\bigskip
Even more speculatively, the foliation pictures suggest a similar
scheme that possibly could eventually yield a good natural proof for the general
geometrization conjecture. The idea is to start with the  unit tangent
bundle, in some Riemannian metric, and turn the fibers into $3$-dimensional
hyperbolic spaces. The aim is to look for a complete flat connection
transverse to the fibers; the affine connection for the Riemannian metric
gives a first approximation, but it is not a complete connection.
In any case, the $\Hy^3$ foliation has a deformation theory parametrized
by quasiconformal structures on its spheres. A relaxation process for
this foliation can be defined, by making the conformal structure
on the sphere at infinity evolve toward the shape of the spheres for its
neighbors in the appropriate direction.  The idea is something like
ripples in a pond when a few equal-sized pebbles are dropped.
The waves spread, perhaps becoming distorted from their initial
conditions if the pond is shallow and irregular, but the wavefronts
stay close to each other.   

Three-manifolds with essential laminations are a special case of this
general pebble-in-the-pond picture: 2-dimensional surfaces in
a 3-manifold allow one to find a connection for
the tangent $\Hy^3$-bundle that is integrable
in two out of three directions, using Candel's uniformization; in this
case, relaxation needs to happen in is only one more direction.

\bigskip
There have been many powerful applications of foliations and laminations
to analyzing 3-manifolds Dehn surgeries on knots.  The phenomena
in this paper suggest that there ought to be a theory which would connect
surgeries along closed orbits of transverse pseudo-Anosov flows to 
rotation numbers and to the Milnor-Wood inequalities.
One way to frame it is this:
most surgeries along closed orbits of a pseudo-Anosov flow
yield manifolds with pseudo-Anosov flows. For which
surgeries is there a transverse foliation?  For which surgeries
does the flow uniformly regulate some transverse foliation?
It seems a reasonable conjecture that every pseudo-Anosov flow
is at least finitely covered by one that admits a transverse foliation.
David Fried (\cite{Fried:AnosovSection}) showed that any
Anosov or pseudo-Anosov flow can be obtained by some Dehn surgery  along
flow lines of the suspension of a pseudo-Anosov homeomorphism of a surface.
To do it, what is required is a section in the complement of some closed
orbits, generalizing the Hedlund-Morse construction.  A question
related to generalizing the Milnor-Wood inequality is to describe
the minimal collections of orbits that need to be removed for the flow
to admit a section.

There are further related questions about tight contact structures:
when is there a contact structure transverse
to a flow, or tangential to a flow?  These questions all seem closely
linked. It would be nice to have a general theory of $\Homeo(\R)$-connections
transverse to flows, or at least, for flows coming from
slitherings of a 3-manifold over a surface; the `nicest' cases are
when the connection has positive or negative curvature, which 
gives a contact structure, or zero curvature, which gives a foliation.
One can similarly look for
analogues of the canonical 1-form in the cotangent bundle, that is,
flows with plane fields (possibly with singularities) that
twist positively, negatively, or not at all.
See \cite{EliashbergThurston} for a discussion of related topics.

Similarly, what happens for surgeries along the leaves
of a foliation (or of an essential lamination?) Is there a
generalization of Fenley's condition to some
class of surgeries along the leaves of a 3-manifold that slithers
around $S^1$? 

\medskip
Given a foliation transverse to the fibers of $M^3 \times S^1$,
is there some finite sheeted covering such that
the pull-back bundle admits a transverse section
(preferably defining a slithering)?  The special
case of the trivial question is the question of whether the three-manifold
virtually fibers over $S^1$. It might be easier to do this when the
foliated bundle is not trivial.  

Is every hyperbolic
three-manifold group isomorphic to a subgroup of $\Homeo \R$?
This is beginning to seem likely.  

\bigskip
In certain mysterious ways,
foliations and essential laminations are quite similar to
hyperbolic structures. Either kind
of structure gives a positive and widely applicable criterion to
show a manifold has many `nice' properties including
infinite fundamental group.  On a surface,
measured laminations can be thought of as the rank one limit of a conformal
structure.  A conformal class of indefinite metrics is a Lorentz cone
structure. The canonical form for a diffeomorphism of a surface in
essence produces one  of the three types of conformal structures invariant
by the diffeomorphism.  Another way to think of the relationship
is that hyperbolic structures are the same thing as
grous acting on the complex $1$-manifold $\CP^1$
that are `taut' in a certain sense.  This is the complex version of
groups that act on $\RP^1$.  Foliations give groups acting on the circle,
with many nice geometric properties.

One can only hope that some day, all these different structures and
constructions will fit together into a single coherent picture.

\bibliography{foliations}
\bibliographystyle{alpha}
\end{document}